\documentclass[10pt]{article}
\usepackage[utf8]{inputenc}
\usepackage{palatino}
\usepackage{comment}
\usepackage{mdframed}
\usepackage{amsmath}
\usepackage{amssymb}
\usepackage{mathrsfs}
\usepackage{amsthm}
\usepackage[toc]{appendix}

\usepackage{graphicx}

\usepackage{hyperref}

\usepackage{todonotes}

\newtheorem{definition}{Definition}
\newtheorem{theorem}{Theorem}
\newtheorem{conjecture}{Conjecture}
\newtheorem{question}{Question}
\newtheorem{remark}{Remark}

\newcommand{\On}{{\rm{Ord}}}
\renewcommand{\emptyset}{\varnothing}

\title{Intrinsic Justification for Large Cardinals and Structural Reflection}

\author{Joan Bagaria,\footnote{Universitat de Barcelona and ICREA. \ttfamily \textbf{bagaria@ub.edu; joan.bagaria@icrea.cat}} \ Claudio Ternullo\footnote{Universitat de Barcelona. \ttfamily \textbf{claudio.ternullo@ub.edu}}}
\date{\today}

\begin{document}

\maketitle

\begin{abstract}

\noindent
We deal with the complex issue of whether large cardinals are
intrinsically justified principles of set theory (we call this the Intrinsicness Issue). In order to do this, we review, in a  systematic fashion, (1.) the abstract principles that have been formulated to motivate them, as well as (2.) their mathematical expressions, and assess the justifiability of both on the grounds of the (iterative) concept of set. A parallel, but closely linked, issue is whether there exist mathematical principles able to yield all known large cardinals (we call this the Universality Issue), and we also test principles for their responses to this issue. Finally, we discuss the first author's Structural Reflection Principles (SRPs), and their response to Intrinsicness and Universality. We conclude the paper with some considerations on the global justifiability of SRPs, and on alternative construals of the concept of set also potentially able to intrinsically justify large cardinals.  
 
\end{abstract}

\section{Introduction}\label{Intro}

Among new set-theoretic axioms, Large Cardinal Axioms hold a place of honour. These axioms have revealed to be immensely successful, both in terms of mathematical (set-theoretic) consequences, and in terms of connections they have helped establish among different areas of mostly, but not exclusively, set-theoretic interest. Proof-theoretically, their investigation has brought in a far-reaching method to calibrate the strength of statements that are independent from the $\mathsf{ZFC}$ axioms.\footnote{For fundamental preliminaries on large cardinals, cf. \cite{kanamori2009}, Introduction.} The study of large cardinals seems to be, in addition, inexhaustible, as new hypotheses, whose usefulness may not, as yet, be fully understood, but will presumably be revealed at some point, keep cropping up in set-theoretic work. 

Their status and justification as new axioms of set theory is, however, a different story, as they do not seem to have a straightforward intuitive appeal. Indeed, \textit{prima facie}, many large-cardinal notions would not seem to follow from the \textit{concept of set}, and, thus, would not be, according to a generally accepted classification of forms of justification for set-theoretic axioms, \textit{intrinsically} justified.\footnote{For the origin of the classification, and its meaning, see \cite{godel1947} and \cite{godel1964}. G\"odel's ideas have been explained in more detail in \cite{wang1996}. Further useful clarifications may be found both in \cite{maddy1996} and \cite{koellner2009b}.} This is precisely what we, for the sake of brevity, shall call:

\vspace{11pt}

\noindent
\textbf{Intrinsicness Issue}. Are large cardinals justified in light of the concept of set?

\vspace{11pt}

\noindent
Now, over the years, several abstract (that is, \textit{intuitive}) motivating principles have been put forward to justify the introduction of large cardinals. Therefore, providing an answer to the Intrinsicness Issue, ultimately, seems to consist in tackling two fundamental questions: (i) \textit{whether} (and \textit{what}) abstract principles used to introduce large cardinals really are \textit{intrinsically justifiable}, and (ii) \textit{how}, that is, on the grounds of \textit{what} concept of set. (i) has, in turn, two parts: (i.a), the justifiability of the abstract motivating principles themselves and (i.b), the justifiability of the mathematical principles expressing them.

This undertaking is very complicated at different levels. As far as (i) is concerned, some authors have taken motivating principles to play a purely \textit{heuristic} (or \textit{explanatory}) role.\footnote{Cf. \cite{maddy1997}, but also \cite{parsons2008}, discussed later on, in section \ref{equilibrium}.} As regards (ii), it has been argued that, if one takes the concept of set to be the \textit{iterative concept of set}, then large cardinals will not be intrinsically justified. An argument to this effect, made by Ralf Schindler in \cite{schindler1994}, may be summarised as follows:

\vspace{11pt}

\noindent
\textbf{Schindler's Challenge}. Believing in the reality (truth) of large cardinals is not justified, since:
\begin{enumerate}
    \item Our only apparently good arguments for the actual existence of large cardinals in the set-theoretic universe, $V$, come from reflection principles.
    \item There are no impredicative\footnote{
    More specifically, Schindler claims that `We cannot believe in the existence of nonpredicative classes on philosophical grounds', cf. \cite{schindler1994}, p. 463. For the distinction between predicative and impredicative classes, see section \ref{absolute}.}  classes.
    \item ${\rm RP}\Pi^1_1$ presupposes the existence of impredicative classes,\footnote{\cite{schindler1994} shows, in $\mathsf{NBG}$, that ${\rm RP}\Pi^1_1$ fails if one assumes that the class-quantifiers range  only over predicative classes. Note that, for the sake of notational consistency throughout the paper, we use ${\rm RP}\Pi^1_1$, instead of \cite{schindler1994}'s ${\rm PR}\Pi^1_1$, to denote Bernays' Reflection.} and
    \item ${\rm RP}\Pi^1_1$ is the weakest reflection schema yielding large cardinals\footnote{${\rm RP}\Pi^1_1$ yields the existence of Mahlo cardinals, but not of weakly-compact cardinals.}
\end{enumerate}
where ${\rm RP}\Pi^1_1$ is Bernays' Reflection Principle, namely the schema asserting that true (in $V$) $\Pi^1_1$ formulae are true in some transitive set.\footnote{\cite{schindler1994}, pp. 458-9. We shall discuss this principle in more detail in Section \ref{Refl}. Schindler's Challenge is also addressed by \cite{mccallum2021}, pp. 199-200, to which we are indebted for the nomenclature. For the full details of \cite{schindler1994}'s argument, see, again, section \ref{Refl}, in particular fn. \ref{schindler}.} Schindler's conclusion is that:

\begin{quote}

[...] anyone who does not believe in nonpredicative classes on philosophical grounds has no justification at all for believing in the existence of large cardinals. (\cite{schindler1994}, p. 459)
    
\end{quote} 

Now, we believe that, if one really wants to make progress on the Intrinsicness Issue, one is bound to address, among other things, Schindler's Challenge and provide a response to it. The philosophical undertaking has already been sketched in \cite{bagaria2023}; in this paper, we wish to elaborate more substantially on the argumentative strategy presented there, and also provide a more comprehensive account of the subtleties, both mathematical and philosophical, involved in the issue.

It is fundamental, already at this stage, to introduce one further issue, which is closely related to, and, we think, also bears on, the potential responses to the Intrinsicness Issue, namely the:

\vspace{11pt}

\noindent
\textbf{Universality Issue}. Is there a set-theoretic principle able to yield \textit{all} known large-cardinal notions?

\vspace{11pt}

\noindent
The connection between the Intrinsicness and the Universality Issues may be explained as follows. On the one hand, one could take \textit{universality} to be a measure of the \textit{intrinsicness} of a motivating principle.\footnote{Ideally, one would like to formally express that a principle $P$ is more \textit{universal} than a principal $Q$, if $P$ \textit{yields} \textit{more} large cardinal hypotheses than $Q$, but the notion of `yielding' is far from being straightforwardly formalisable. However, later in the text, we shall encounter some   examples of principles enjoying `greater universality' than others.} On the other hand, a principle which met all conditions for being seen as \textit{intrinsically justified}, but which had a low degree of \textit{universality} would hardly fulfill general justificatory purposes with respect to large cardinals.

As a consequence, overall, the goal of this paper is to investigate the Intrinsicness Issue in connection with fundamental abstract motivating principles for large cardinals, with a view to successfully meeting Schindler's Challenge, but also to responding to the Universality Issue. In the end, we will make the case that Structural Reflection Principles, introduced and investigated by the first author, yield an answer to Schindler's Challenge, are intrinsically justified (if `strong reflection' principles are, in general, thus justified), and meet, to a very large extent, the Universality Issue. 

\medskip

Here's a short summary of the paper's contents. In section \ref{sec2}, we review the concept of set, and the issue of the necessary use of classes to motivate large cardinals. In section \ref{sec3}, we examine \textit{three} fundamental abstract principles for large cardinals, and show that they exhibit an increasing degree of universality that reaches its apex with Resemblance. In section \ref{SRP}, we discuss Structural Reflection Principles. Finally, in section \ref{secfinal}, we show how large cardinals may turn out to be intrinsically (and globally) justified under Structural Reflection (section \ref{equilibrium}). In passing, in \ref{interlude}, we also briefly consider non-conservative types of responses to the Intrinsicness Issue, that is, responses hinged on an alternative construal of the concept of set. 

\section{The Iterative Concept of Set, Classes and Large Cardinals}\label{sec2}

In this section, we begin our investigation of the question of whether large cardinals are justified by the (iterative) concept of set; in particular, we show that the concept of set alone does not seem sufficient to motivate any of them. This is because even the weakest among the known large-cardinal notions seem to require the grasp of notions, such as that of Cantor's Absolute, which \textit{transcend} the scope of the (iterative) concept of set. 

The discussion to be found here may also be taken to set the stage for the more substantial discussion of abstract motivating principles, Structural Reflection Principles, and the intrinsic justification of large cardinals through these, in sections \ref{sec3}, \ref{SRP} and \ref{secfinal}. 

\subsection{The Iterative Concept}\label{IC}

It is very often repeated, and held, both by practising set-theorists and philosophers of mathematics that `concept of set' means the:

\vspace{11pt}

\noindent
\textbf{Iterative Concept of Set (IC)}. Sets are formed in \textit{stages} within an iterative, and cumulative, process indexed by the ordinal numbers.

\vspace{11pt}

\noindent
By IC, the study of sets  reduces to the study of the members of the well-founded hierarchy $V$, indexed by the ordinal numbers. To review, one defines $V_0=\emptyset$, and then takes power-sets at \textit{successor levels}, that is, one defines: $V_{\alpha +1}=\mathcal{P}(V_\alpha)$, and the union of all previous levels at \textit{limit levels}, i.e., if $\lambda$ is a limit ordinal, then  $V_{\lambda}= \bigcup_{\alpha<\lambda} V_\alpha $. In $\mathsf{ZF}$,  via the Axiom of Foundation, then one proves that every set $x$ belongs to some $V_\alpha$.\footnote{Although, arguably, IC is foreshadowed in \cite{godel1944} and \cite{godel1947}, the fundamental traits of the conception are described in \cite{boolos1971} (see also \cite{boolos1989}) and \cite{wang1974}, and further discussed in \cite{parsons1977} and \cite{shoenfield1977}. For a more recent overview, see \cite{incurvati2020}, Ch. 2-3.} 

All of the $\mathsf{ZFC}$ axioms have been taken to be true of this concept.\footnote{Cf. \cite{wang1974}, p. 184, and \cite{boolos1971}, p. 499ff.} The Axiom of Replacement, though, is a controversial case: the `regressive' argument according to which without Replacement one would thwart the development of set theory, since, then, already $V_{\omega+\omega}$ couldn't be proved to exist, has been taken by some authors to be insufficient.\footnote{Cf. \cite{potter2004}, in particular pp. 218ff. The two main sources for IC, \cite{boolos1971} and \cite{wang1974}, differ on the assessment of the justifiability of Replacement via IC. See also \cite{incurvati2020}, section 3.6.} For another justificatory strategy, the Reflection Principle, with which we will be concerned in a moment, has also been seen as justifying Replacement.\footnote{Cf. \cite{incurvati2020}, pp. 95-100, and \cite{potter2004}, pp. 221ff.} But then, again, since Reflection, as we shall see, does not straightforwardly fall within the compass of IC, there are further reasons to doubt that IC is able to justify (relevant) bits of $\mathsf{ZFC}$.

But for the sake of our discourse, let us assume that IC is, indeed, able to justify all of the $\mathsf{ZFC}$ axioms. Now, can IC also justify the introduction of large cardinals? In order to answer this question let us consider the large-cardinal notion of strongly  inaccessible cardinal:

\begin{definition}[Strongly Inaccessible Cardinal]

A cardinal $\kappa$ is strongly inaccessible if it is: uncountable, strong limit (i.e., for all cardinals  $\lambda<\kappa$, also $2^{\lambda}<\kappa$), and regular (i.e., $cf(\kappa)=\kappa$). 

\end{definition}

\noindent
along with the simplest abstract motivating principle which seems to be able to motivate it, i.e. Reflection, expressed, mathematically, as follows::
\vspace{11pt}

\noindent
\textbf{Reflection Principle}. If $V$ has some property $P$, then there exists an ordinal $\alpha$ such that $V_{\alpha}$ has $P$.\footnote{For an introductory, but more rigorous treatment of the notion of `$V$ having a property $P$' (\textit{vis-à-vis} `$V$ satisfying some sentence $\varphi$'), see \cite{bagaria2023}, section 1.1.} 

\vspace{11pt}

\noindent
 Then one reasons as follows. In $V$, the class $\On$ of all ordinal numbers is uncountable, closed under cardinal exponentiation, and `regular', i.e.,  there is no  definable (with parameters) sequence of ordinals, indexed by an ordinal, that is cofinal in $\On$ (this just follows from Replacement).  Then by the Reflection  Principle there is some $V_\kappa$, with $V_\kappa =H_\kappa$, $\kappa$ an uncountable cardinal, which reflects these properties, i.e., they hold in $V_\kappa$. In particular, $\kappa$ is uncountable and closed under cardinal exponentiation, but not necessarily regular.  However, if we take the property: `$\On$ is regular' to be a second-order property of $V$, and this second-order property holds in $V_\kappa$, then, since \emph{all} ordinal sequences indexed by an ordinal less than $\kappa$ into $\kappa$,  are in $V_{\kappa +1}$,  $\kappa$ must indeed be regular,  and  so \textit{strongly inaccessible}. Through iterating this strategy, one may also introduce hyper-inaccessible cardinals, hyper-hyper-inaccessible cardinals, and so on, up to the level of \textit{Mahlo cardinals}.\footnote{A cardinal $\kappa$ is \textit{Mahlo} if and only if it is regular and the set of strongly inaccessible cardinals below $\kappa$ is \textit{stationary}. A subset $A$ of  $\kappa$  is called \textit{stationary} if it intersects all   club (closed and unbounded) subsets of $\kappa$.} 

We have thus reduced the issue of the intrinsic justifiability of strongly inaccessible cardinals to the issue of the intrinsic justifiability of the Reflection Principle, construed as: (1) allowing second-order properties of $V$ 
to be reflected to some $V_\alpha$, and (2) having second-order quantifiers ranging over \emph{all} subsets of $V_\alpha$, not just the \textit{definable} ones. Now the question is: is IC able to justify the Reflection  Principle, stated in this form? It wouldn't seem so, as the set-theoretic universe envisaged by IC is never a \textit{completed totality} (as otherwise it would form a new stage of the iterative process), and so neither are the \textit{subclasses} of $V$. Hence, it would seem that, by IC's own lights, second-order properties of $V$ can only be taken to be legitimate when second-order quantification ranges only over \textit{definable subclasses} of $V$, but this, as we have seen, is insufficient to intrinsically justify even strongly inaccessible cardinals.\footnote{See also \cite{bagaria2023}, section 1.}

If a bit crudely, this example already showcases the difficulties one encounters in attempting to justify large cardinals intrinsically, that is, exclusively based on IC.

At this point, some readers may object that there exists a different, but maybe more promising, way to re-interpret (and re-formulate) the  Reflection Principle based on set-theoretic \textit{potentialism}, that is, the idea that quantification over \textit{all} sets is an \textit{indeterminate} notion (precisely because there is no \textit{completed} totality of all sets). In particular, versions of potentialism assert that $V$ is an `unfinished' object, that the variables in the axioms of set theory do not range over \textit{all} of $V$, but just over some fixed initial segment of it, that what sets there are is not known \textit{a priori}, or that there \textit{may} be more sets than one finds in some provisionally fixed universe of sets.\footnote{For a primer on the actualism/potentialism divide, it is useful to read \cite{koellner2009b}'s section 1.} We defer discussion of potentialist conceptions, and of the difficulties connected to them, to sections \ref{Refl} and \ref{RR}, when it will be fully clear how the Reflection Principle could be formulated mathematically.

\subsection{The Absolute Infinite and Classes}\label{absolute}

The historical details of how other collections beyond sets gradually manifested themselves in the set-theoretic landscape are well-known, but, for our purposes, it is worth briefly recalling them. 

Presumably as a consequence of the discovery of the paradoxes, Cantor formulated the doctrine that there existed `multiplicities' the `being together' of whose elements was \textit{inconsistent}: such collections, of `all ordinals', `all sets', etc. he called \textit{inconsistent}, or \textit{absolutely infinite} multiplicities, as opposed to `consistent multiplicities', that is, sets.\footnote{Cantor's conception of the Absolute is first presented in \cite{cantor1883c}, then more sharply outlined in \cite{cantor1885b} (as well as in \cite{cantor1887}, \cite{cantor1888}). For the notion of `inconsistent multiplicity', see Cantor's 1899 letter to Dedekind in \cite{ewald1996b}, pp. 930-5. An exhaustive discussion of the conception is in \cite{hallett1984} and \cite{jane1995}; for further historical details, see \cite{ferreiros2004} and \cite{tapp2012}.} Now, although the overall interpretation of (the value of) the conception is controversial, Cantor is often credited with having also held that:   

\begin{enumerate}
    \item Sets may not be sufficient to investigate, and fruitfully develop, the whole of set theory
    \item Absolutely infinite entities are \textit{determinately existing} objects as much as sets\footnote{Arguably, there is an unresolved tension between Cantor's claiming that absolutely infinite collections are \textit{inconsistent}, and the parallel claim that they conceivably \textit{exist}. On this issue, cf. \cite{jane1995}.}
    \item The Absolute cannot be \textit{measured} (for that matter, cannot even be fully \textit{described})
    
\end{enumerate}

\noindent
Several scholars have advocated the view that Cantor's conception was, indeed, mathematically fruitful, as it led to the formulation of what came to be known as:

\vspace{11pt}

\noindent
\textbf{Limitation of Size Conception (LSC)}. A multiplicity $M$ is a set if and only if it isn't `too big', that is, if and only if $M$ is \textit{smaller than} $V$.

\vspace{11pt}

\noindent
the basis of, respectively, \cite{vonneumann1925}'s set/class theory (that explicitly contains what may be called an Axiom of Limitation of Size),\footnote{This is the Axiom IV.2 of \cite{vonneumann1925}, p. 400. Note that LSC above is a restatement of that axiom. A thorough discussion of von Neumann's theory, also of its historical aspects, may be found in \cite{fraenkel1973}, Ch. II.7.}  but also, arguably, of \cite{zermelo1908a}'s early axiomatic theory $\mathsf{Z}$, in particular, of the Axiom of Separation.\footnote{The point is made in, among other works, \cite{fraenkel1973}, p. 32, and criticised in \cite{hallett1984}, pp. 198ff.} 

\medskip

Among the subsequent elaborations of Cantor's doctrine of the Absolute liable to be linked to large cardinals, one may consider the following principle, a restatement of bullet point (3.) above:

\vspace{11pt}

\noindent
\textbf{(Unknowability)}. The Absolute is unknowable. 

\vspace{11pt}

\noindent
As we shall see in more detail in section \ref{sec3}, according to Wang, G\"odel expressed the view that all the axioms of set theory should be reduced to Unknowability (equivalently, to the Reflection Principle, taken to be Unknowability's mathematical embodiment).\footnote{Cf. \cite{wang1996}, p. 283. See later section \ref{RR}.}  In a similar fashion, \cite{reinhardt1974b}, expanding on work in \cite{ackermann1956}, takes Cantor's Absolute and, to some extent, Unknowability, as the starting point of his radical extension of the Reflection Principle. Finally, more recent work on large cardinals, as we shall see in a moment, is still strongly indebted to Unknowability and Cantor's doctrine of the Absolute.

\medskip

The progression from set to class theory may also be justified in the way sketched by Parsons below:

\begin{quote}
    This process would be that of gradually imposing on our discourse an interpretation which makes the original universe a set. Even before the introduction of classes, the application of classical logic statements about all sets could be taken as a first step in this direction. (\cite{parsons1974}, p. 219)
\end{quote}

\noindent
that points to the essential role of classical logic and is more amenable to \textit{potentialist} ideas one may already find enunciated in \cite{zermelo1930}.\footnote{\label{zermelo}\cite{zermelo1930} holds that the set/class distinction is a \textit{temporary} one, in the sequence of normal domains (models of $\mathsf{ZFC}_2$) $D= V_\kappa$, $D'=V_{\kappa+1}$, $D''=V_{\kappa+2}$, ..., where $\kappa$ is the least \textit{strongly inaccessible} cardinal: in particular, for all ordinals $\alpha$, the \textit{classes} of $D^{\alpha}$ are \textit{sets} in the next domain $D^{\alpha+1}$.}   

In any case, granting (at least, provisionally) that the `addition' of classes is justified, either on account of Cantor's, or of Parsons' (Zermelo's) conception, now the choice is between a \textit{predicative} and an \textit{impredicative} theory of classes, according to whether one allows for the introduction of a Predicative or Impredicative Class Comprehension Principle. Class Comprehension Principles are axioms of this sort:

\[ (\exists X)(\forall y)(y \in X \leftrightarrow \Phi) \tag{Comp} \]

\noindent
where the second-order variables (predicates) range over classes, and the first-order variables over sets. If $\Phi$ does not contain class quantifiers, then the resulting theory ($\mathsf{NBG}$) is \textit{predicative}, if it does have class quantifiers, then the theory ($\mathsf{MK}$) is \textit{impredicative}. 

For our purposes in the next sections, it is fundamental to notice that Predicative Class Theory is proof-theoretically equivalent to $\mathsf{ZFC}$; in fact, one may just view it as a more convenient way to express the \textit{schematic} reference to properties in the context of $\mathsf{ZFC}$. The ensuing attitude about classes, which takes these to be just metatheoretic formulas, has been called \textit{definabilism}.\footnote{See later in the text, section \ref{interlude}.} In contrast, impredicative class theory seems to presuppose that there exists a determinate, self-standing, realm of entities alongside sets, \textit{the} classes, and some authors have seen this conception as befitting, even originating with, Cantor's doctrine of the absolute infinite.\footnote{This is, for instance, \cite{horsten-welch2016}'s standpoint, which we shall briefly review in section \ref{RR}.}

Whatever the choice at this stage, what is clear is that, if one is willing to consider the Absolute as an integral part of the concept of set, then one has it that the  Reflection Principle, in its strong form, may be justified by this (enriched) concept of set and, as a consequence, that also strongly inaccessible cardinals in the examples above, and presumably some other types of large cardinals, may be \textit{intrinsically justified}. But, so far, no cogent reason to take this step has emerged; later, in section \ref{interlude}, we will briefly discuss ways in which this approach may be put into effect.

\section{The Motivating Principles}\label{sec3}

In this section, we carry out the bulk of the tasks we have set ourselves in the Introduction relating to questions (i.a), (i.b) and (ii), that is, to the intrinsic justifiability of abstract motivating principles and of the corresponding mathematical principles. For each motivating principle, we will also indicate how well (or badly) it fares with the Universality Issue. 

As far as the latter is concerned, we will also show that there has been a clear evolutionary trend in the formulation of motivating principles, reaching a maximum with Resemblance (section \ref{resemblance}).

\subsection{Reflection}\label{Refl}

We have seen (in section \ref{IC}) that the Reflection Principle, in the strong form that allows for the reflection of second-order properties, may not be licensed by IC, insofar as it seems to be essentially transcending IC's conceptual resources.
 In what follows, we wish to look into the issue more closely, by examining the concrete mathematical principles motivated by the Reflection Principle, and their consequences in terms of large-cardinal strength.

\medskip

An essential feature of the set-theoretic universe, $V$, is that it reflects any statement holding in it to some of its rank-initial segments. This fact is a theorem of $\mathsf{ZF}$ (\cite{levy1960}), the Reflection Theorem, namely the following schema of statements, one for each formula $\phi(x_1\ldots ,x_n)$ of the language of set theory:\footnote{In fact, a stronger form of Reflection is provable in $\mathsf{ZF}$. Namely, for each  $n$ there exists a closed and unbounded proper class $C^{(n)}$ of ordinals such that $V_\alpha$ is a $\Sigma_n$-elementary substructure of $V$, for every $\alpha$ in $C^{(n)}$. It should be noted that, provably in $\mathsf{ZF}$, the RP will hold in any cumulative hierarchy, i.e., in any class of sets indexed by a club class of ordinals which forms a chain under inclusion and is continuous (i.e., unions are taken at limit points).}

\[ \exists \alpha \forall x_1, ..., x_n \in V_\alpha (\phi(x_1, ..., x_n) \leftrightarrow \phi^{V^{\alpha}}(x_1, ..., x_n)) \tag{RP$_0$}\]

Now, it is known that strengthenings of RP$_0$ are proof-theoretically equivalent to \textit{some} large cardinal notions. For instance, if we allow $\phi(x_1,...,x_n)$ to contain second-order parameters then we get an RP which holds in $V_\kappa$  if and only if $\kappa$ is an \textit{inaccessible} cardinal. If, in addition, we require that $\alpha$ is \textit{inaccessible}, then we get:

\[\exists \alpha (Inac(\alpha) \wedge \forall x_1, ..., x_n \in V_\alpha (\phi(x_1, ..., x_n) \leftrightarrow \phi^{V^{\alpha}}(x_1, ..., x_n))) \tag{RP$_1$}\]

\noindent
(where \textit{Inac}$(\alpha)$ = `$\alpha$ is inaccessible') which holds in $V_{\kappa}$ if and only if $\kappa$ is \textit{Mahlo}. By iterating this strategy, one obtains RP's equivalent to such cardinals as \textit{Mahlo}, \textit{weakly compact} and $\Pi^m_{n}$ \textit{indescribable cardinals}.\footnote{For the definitions of weakly compact and $\Pi^{m}_{n}$-indescribable cardinals, see \cite{kanamori2009}, respectively, pp. 37ff. and 56ff.}

All of these large-cardinal notions are globally motivated by the following second-order principle formulated in \cite{bernays1976} (and, as a consequence, also known as Bernays' Reflection Principle): 

\[ \Phi \rightarrow (\exists u)(Trans(u) \wedge \Phi^{u}) \tag{${\rm RP}\Pi^1_1$} \]

\noindent
for all $\Pi^{1}_{1}$-formulae $\Phi$, where $u$ is a set variable and $Trans(u)$ means: `$u$ is transitive'. As readers already know, this is the principle  considered by Schindler in \cite{schindler1994} (see Section \ref{Intro}).

Now, as first noticed by Reinhardt, if one allows for  third-order parameters in the formula $\Phi$, the resulting RP is inconsistent. More specifically, as shown by Koellner, every such strengthening of ${\rm RP}\Pi^1_1$, if consistent, cannot yield any large cardinal up to or beyond  the first $\omega$-Erd\H{o}s cardinal, $\kappa (\omega)$.\footnote{\cite{koellner2009}, p. 210ff.} Koellner comments on this result as follows:

\begin{quote}

..our main limitative result is also schematic and the proof would appear to be able to track any degree of reflecting on reflection, the Erd\H{o}s cardinal $\kappa(\omega)$ appears to be an impassable barrier as far as reflection is concerned. This is not a precise statement. But it leads to the following challenge: Formulate a strong reflection principle which is intrinsically justified on the iterative conception of set and which breaks the $\kappa(\omega)$ barrier (\cite{koellner2009}, p. 217.).\footnote{For the definition of Erd\H{o}s cardinals, see \cite{kanamori2009}, p. 80.} 
    
\end{quote}

\medskip

\noindent
Now, it would seem that what we may call Koellner's Challenge cannot be successfully met by any Reflection Principle, and the reason, as we have seen, is the one pointed out by \cite{schindler1994}: ${\rm RP}\Pi^1_1$ already implies the existence of \textit{impredicative classes} that are not licensed by IC.\footnote{\label{schindler}In particular, \cite{schindler1994}, pp. 460-1, shows that ${\rm RP}\Pi^1_1$ proves $\Delta^{1}_{1}$-Comprehension in $\mathsf{NBG}$, and the latter fails if the domain consists only of predicative classes.} \textit{A fortiori}, any stronger, and consistent, reflection principle is bound not to be justifiable solely in the light of IC.  

But even if large cardinals really happened to be justifiable exclusively in the light of IC, Koellner's result shows that Bernays' Reflection may not be improved by any other principle potentially able to capture \textit{all} of them, so also the other part of Koellner's Challenge cannot be met by this kind of Reflection.

\medskip

At the end of section \ref{IC} we have hinted at the possibility that a different, \textit{potentialist}, construal of $V$ may successfully overcome the difficulties relating to Reflection. The strategy has been laid out by William Tait. In particular, \cite{tait2005b} presents a strategy to introduce small large cardinals, namely those large cardinals compatible with the axiom $V=L$, in a \textit{bottom-up} fashion, in a way, that is, which does not presuppose that $V$ is a \textit{completed} totality.  

On Tait's view, Cantor's doctrine of the Absolute is obscure (and Cantor's set/class distinction, which is based on it, question-begging),\footnote{Cf. \cite{tait2005b}, p. 133.} so, rather than using, or referring to, absolutely infinite collections, one should just use \textit{initial segments} thereof, in particular initial segments of $\On$, whose existence is, instead, licensed by IC, and view large cardinals as arising from existence conditions  $C$ (expressed by specific formulas) associated to those initial segments (in practice, ordinals). For appropriately selected $C$'s, we will have that an ordinal $\alpha$ is inaccessible, Mahlo, $\Pi^{m}_{n}$-indescribable, etc.

Based on \cite{tait2005b}'s approach, \cite{koellner2009b} shows that existence conditions are precisely exhausted by a restricted class of second-order formulas with higher-order parameters, the  formulas $\Gamma^{(2)}$ considered by Tait; Tait's bottom-up approach would, thus, show that, at least, the large cardinals obtained by reflection of $\Gamma^{(2)}$-formulas are justified by IC, so Schindler's Challenge is met. 

Recently, \cite{mccallum2021} has considered an extension of Tait's strategy that may even be able to justify stronger large-cardinal notions. In particular, McCallum has suggested that \cite{roberts2017}'s theory $\mathsf{ZFC2}_{S}+R_{S}$ (that is, second-order $\mathsf{ZFC}$ with a satisfaction predicate $S$ and reflection principle $R_{S}$ for all formulas in the language of $\mathsf{ZFC2}_{S}$), if suitably strengthened to an $\omega$-th order theory $\mathsf{ZFC}_{\omega}+R_{\omega}$ ($R_{\omega}$ being the $\omega$-th order generalisation of $R_{S}$) is able to yield all $n$-extendible cardinals, with $n$ any positive integer; in turn, the latter, is, under the stipulation that only certain `reflecting structures' are taken into account, equivalent to an $\omega$-th order theory, $T$, that McCallum introduces to incorporate all of Tait's $\Gamma^{(2)}$ formulas.\footnote{\cite{mccallum2021}, Theorem 3.3. Also see the discussion on p. 209ff.} Given the established analogy between $T$ and $\mathsf{ZFC}_{\omega}+R_{\omega}$, and on the grounds of the fact that the latter yields all $n$-extendible cardinals, for any integer $n$, McCallum argues that Tait's strategy may be invoked to intrinsically justify also \textit{extendible} cardinals. So, if McCallum's argument is sound, then we would have that Tait's bottom-up approach, which, Tait argues, is exclusively based on IC, would vindicate the justifiability even of extendible cardinals, hence Koellner's (and, of course, also Schindler's) Challenge would be met. 

We have several qualms about the tenability of both Tait's and McCallum's strategies. 

On the one hand, as already pointed out by Koellner, there are problems of (i) excessive generality and (ii) consistency with regard to Tait's approach, since, in particular, some existence conditions give rise to inconsistent reflection principles.\footnote{\cite{koellner2009b}, p. 208.}  

On the other hand, the restriction of Reflection to particular formulas, i.e. Tait's $\Gamma^{(2)}$ formulas, and McCallum's focus on theories and principles, such as \cite{roberts2017}'s, that incorporate \textit{satisfaction predicates}, finally, McCallum's own stipulation that certain `reflecting structures' be ruled out, in order for $T$ and $\mathsf{ZFC}_{\omega}+R_{\omega}$ to be really equivalent,\footnote{\cite{mccallum2021}, p. 209.} do not seem to be licensed by IC. 
As a consequence, even if the bottom-up approach successfully prevented the difficulties discussed in section \ref{IC}, and even if it really were a valid alternative to the `actualist' view of $V$ discussed in section \ref{absolute},  it would still be questionable that the large cardinals justified by the corresponding reflection principles really are justified by IC. 

\subsection{Uniformity}\label{uniformity}

Uniformity is one of G\"odel's five criteria for the introduction of new axioms in \cite{wang1974} and \cite{wang1996}. This principle is also discussed in, among other works, \cite{solovay-reinhardt-kanamori1978} and \cite{kanamori-magidor1978} and, finally, \cite{maddy1988a}. All of these works take Uniformity to correspond to a process of `generalisation of the properties of sets to other sets'. The quintessential expression of this is in \cite{solovay-reinhardt-kanamori1978}:

\begin{quote}
    
[I]t is in many ways quite reasonable to attribute
certain properties of $\omega$ to uncountable cardinals as well, and these considerations can yield the measurable and strongly compact cardinals. Also, in considerations involving measurable cardinals, natural strengthenings of closure properties on
ultrapowers yield the supercompact cardinals [...] (p. 75)\footnote{Cf. also \cite{kanamori-magidor1978}, p. 104. \cite{maddy1988a} discusses a kindred principle called `whimsical identity', the idea that, if, for instance, $\omega$ were the only cardinal to have the properties it has, then $\omega$ would be whimsical. Cf. \cite{maddy1988a}, p. 502.}
\end{quote}

Here the authors refer to $\omega$'s being, for instance, \textit{inaccessible} (since it is regular and limit) and also to its being \textit{measurable} (since there exists an $\omega$-complete non-principal ultrafilter on it), as the basis for the existence of inaccessible and measurable cardinals, the only difference with $\omega$ being that they are required to be uncountable. Another, less frequently discussed, example made by the authors is that of:

\vspace{11pt}

\textbf{Supercompact Cardinals}. $\kappa$ is supercompact if it is $\gamma$-supercompact for all ordinals $\gamma$.

\vspace{11pt}

\noindent
which generalise over:

\vspace{11pt}

\noindent
\textbf{$\gamma$-Supercompact Cardinals}. $\kappa$ is $\gamma$-supercompact iff there exists an elementary embedding $j: V \rightarrow M$ such that $\kappa$ is the critical point of $j$, $j(\kappa)$ is greater than $\gamma$, and  $M$ is closed under $\gamma$-sequences. 

\vspace{11pt}

\noindent
and, in turn, generalise over elementary embeddings used to define \textit{measurable cardinals} by taking $M$ to be closed under \textit{arbitrary sequences} of length $\gamma$.

Another example, occurring in the context of infinitary languages, is that leading from \textit{weakly} to \textit{strongly compact} cardinals. To review, $\kappa$ is weakly compact if and only if any collection of sentences of the infinitary language $L_{\kappa, \kappa}$, using at most $\kappa$ non-logical symbols, which is $\kappa$-satisfiable, is satisfiable. This easily generalises to $L_{\kappa, \kappa}$, $\kappa$ strongly compact, if the restriction that the collection of sentences use at most $\kappa$ non-logical symbols is lifted. 

For a final, more recent example, Structural Reflection Principles (with which we will extensively deal in section \ref{SRP}), express the transition between \textit{strong} and \textit{Woodin} cardinals, on the one hand, and between \textit{supercompact} cardinals and \textit{Vop\u{e}nka's principle} as the generalisation of, respectively, the principle called $\Pi_1$-SR to $\Pi_n$-SR (for all $n$) and the principle $\Pi_1$-PSR to $\Pi_n$-PSR (for all $n$).\footnote{\cite{bagaria2023}, p. 21ff.} 

At the intuitive level, what would underlie all such generalisations is:

\begin{quote}

[a] process of reasonable induction from familiar situations
to higher orders, with the concomitant confidence in the recurring richness of the cumulative hierarchy.\footnote{\cite{kanamori-magidor1978}, p. 104.}
    
\end{quote}

\noindent
Generalisations are routine in mathematics. There is, thus, no need to see set-theoretic generalisations as having a special status. In particular, it does not seem that the generalisations we have mentioned directly stem from the concept of set, i.e., are licensed by IC. Hence, these particular instances of Uniformity do not make the corresponding large-cardinal notions justified by IC.

However, in the passage above, the `richness of the cumulative hierarchy' is invoked, something which does adumbrate a more substantially motivated use of Uniformity. 

In particular, G\"odel's original enunciation of the principle (the most enigmatic among the five G\"odelian criteria) is broader, and, again, points to a recurring richness of $V$ that foreruns stronger generalisations of Reflection:

\begin{quote}

The universe of sets does not change its character substantially as one goes from smaller to larger sets or cardinals; that is, the same or analogous states of affairs reappear again and again (perhaps in more complicated versions). In some cases, it may be difficult to see what the analogous situations or properties are. But in cases of simple and, in some sense, ``meaningful'' properties it is pretty clear that there is no analogue except the property itself. [...] For axioms of infinity this principle is construed in a broader sense. It may also be called the ``principle of proportionality of the universe": analogues of properties of small cardinals by chance lead to large cardinals. For example, measurable cardinals were introduced in this way. People did not expect them to be large.\footnote{\cite{wang1996}, p. 281.}

\end{quote}

We shall briefly review intuitions related to this second type of Uniformity in subsection \ref{RR}, in particular in the context of Reinhardt's Reflection. Then, it will be reasonable to wonder whether this criterion provides any intrinsic justification for large cardinals. 

\subsection{Resemblance}\label{resemblance}

\subsubsection{Preliminaries. Resemblance as Richness.}\label{VP}

Both Reflection- and Uniformity-inspired considerations very naturally lead to, and culminate in, what may, alternatively, be seen as a strengthening of Reflection or an iterated application of G\"odel's second type of Uniformity, that is, Resemblance.  

This is discussed, along with other principles, in both \cite{solovay-reinhardt-kanamori1978} and \cite{kanamori-magidor1978}. \cite{solovay-reinhardt-kanamori1978} explicates the principle as follows:

\begin{quote}

Because of reflection considerations and, generally speaking, because the cumulative hierarchy is neutrally defined in terms of just the power set and union operations, it is reasonable to suppose that there are $\langle V_\alpha, \in\rangle$'s which resemble each other. 

The next conceptual step is to say that there are elementary embeddings $j:\langle V_\alpha, \in\rangle \rightarrow \langle V_\beta, \in\rangle$. Since this argument can just as well be cast in terms of $\langle V_{f(\alpha)}, \in , X(\alpha)\rangle$'s, where $f(\alpha)$ and $X(\alpha)$ are uniformly definable from $\alpha$, the elementary embeddings may well turn out not to be the identity (p. 75).

\end{quote}

\noindent
So, the key idea underlying the principle is to express the \textit{resemblance} between initial segments of the cumulative hierarchy using \textit{elementary embeddings}. 

Now, one way to justify  Resemblance is through linking it to a slightly different principle, i.e., Richness, expressing, like Uniformity, what Kanamori, in the previous section, referred to as the `recurring richness of the cumulative hierarchy'. The argument can be found in \cite{maddy1988b}, which takes the justifiability of:   

\vspace{11pt}

\noindent
\textbf{Vop\u{e}nka's Principle (VP)}. For any proper class $\mathcal{C}$ of structures of the same type, there exist two different structures $A$ and $B$ in $\mathcal{C}$ that  \textit{resemble} each other, i.e., such that one is \textit{elementarily embeddable} into the other. 

\vspace{11pt}

\noindent
to be a suitable case study. \cite{solovay-reinhardt-kanamori1978} had already taken VP to be motivated by Resemblance, insofar as this aptly reflects the idea, expressed in the quote above, that there might be two different initial segments of $V$ that resemble each other. Now, \cite{maddy1988b} suggests that: 

\begin{quote}

The rule of thumb usually cited as lying behind
this principle [VP, our note] is the idea that the proper class of ordinals is extremely rich [...]. Suppose, for example, that a process is repeated once for each ordinal -- $\On$-many times, we might say -- and every step produces a structure. Then richness implies that no matter how closely we keep track of the
structures generated, there are so many ordinals 
that some will be indistinguishable. (\cite{maddy1988b}, p. 750)\footnote{Incidentally, it should be noted that, for Maddy's suggested procedure to work, the structures produced at every step (or at least at proper-class many steps) need to have the same type. However, as VP is equivalent to VP restricted to structures of the form $\langle V_\alpha , \in, A\rangle$, with $A$ being a constant, this is not a problem.}

\end{quote}

So, the existence of two different but very similar initial segments of $V$, mathematically expressed by the existence of an elementary embedding of one of them into the other, would be made possible by the internal richness of $V$ itself. Now, compare \cite{maddy1988b}'s statement of Richness to \cite{martin1976}'s articulation of Reflection:

\begin{quote}

Reflection principles are based on the idea that the
class $On$ of ordinal numbers is so large that, for any reasonable property $P$ of the universe of all sets $R_{On}$, $On$ is not the first stage such that $R_{\alpha}$ has $P$. (\cite{martin1976}, p. 85)

\end{quote}

Again, here Martin  refers to the proper class of ordinals, $\On$, as being extremely `rich', but views this characteristic of $\On$ as being, in fact, expressible in terms of Reflection.  Thriving on this, \cite{maddy1988b} explains that one may equivalently use Reflection and Richness in order to get Resemblance and, as a consequence, justify VP. 

One may think that \cite{maddy1988b}'s argument, based on Richness, might help us make the case of the \textit{intrinsic justifiability} of Resemblance straight away; however,  as is clear from looking at \cite{martin1976}'s formulation of Reflection, and from Maddy's argument itself, in both there is an explicit  reference to the proper class of all ordinals, $\On$, which seems to presuppose conceptual resources beyond those of IC. 

By contrast, \cite{solovay-reinhardt-kanamori1978}'s reference to mutually resembling $V_{\alpha}$'s does not directly require of one to refer to $V$, $\On$, or other proper classes, and, thus, in this form, Resemblance may be taken to be compatible with IC. In this case, however, what would be left to be established is that the existence of elementary embeddings between different initial fragments of $V$ really is justified by IC.

\subsubsection{Elementary Embeddings as Self-Similarity Principles}\label{self}

A noticeable intuition motivating Resemblance, in fact, the use of elementary embeddings, is what we shall call the `self-similarity of $V$'; thus, an accurate analysis of this intuition will gradually lead us to envisage elementary embeddings of a distinct kind as the correct way to flesh out Resemblance. The self-similarity of $V$ may be expressed informally as follows: 

\vspace{11pt}

\noindent
\textbf{(Self)}. The universe of sets contains `copies' of itself.  

\vspace{11pt}

If by `copies' we just mean any isomorphic copies, then Self is true in a very strong  sense. Indeed, for any set $a$, the class 
$$V\restriction a:=\{ x: a\in TC(\{ x\})\}$$
is isomorphic to $V$, via the map $i_a:V\to V\restriction a$ defined recursively (with $a$ as a parameter) by:
\begin{enumerate}
    \item[] $i_a(\emptyset)=a$
    \item[] $i_a(x)=\{ i_a(y):y\in x\}$.
    \end{enumerate}
Observe that the class $V\restriction a$ is both $\Sigma_1$ and $\Pi_1$ definable, with $a$ as a parameter, and the map $i_a$ is $\Sigma_1$-definable, also with $a$ as a parameter.

Note however that $V$ cannot be isomorphic to a \emph{transitive} class different from $V$, for, provably in $\mathsf{ZF}$, any two transitive isomorphic classes must be  identical. Thus, adding the requirement that the copies need to be transitive makes Self vacuous, for the only transitive isomorphic copy of $V$ contained in $V$ is $V$ itself. 

But then, a natural weaker requirement would be that there is an isomorphic copy, $V'$, of $V$, different from $V$ and hence not transitive, which is very \emph{similar} to a transitive class $M$. A natural way to construe this notion of similarity in a very precise manner is to require that $V'$ is  an \textit{elementary substructure} of $M$. Now note that if $\pi:V\to V'$ is an isomorphism and $V'$ is an elementary substructure of $M$, then ${\rm Id}\circ \pi :V\to M$ is an elementary embedding.\footnote{To review, an elementary \emph{embedding} $j: \mathfrak{A} \rightarrow \mathfrak{B}$ of a structure $\mathfrak{A}$ into a structure $\mathfrak{B}$, both of the same type, is a function from the domain $A$ of $\mathfrak{A}$ into the domain $B$ of $\mathfrak{B}$ such that for every formula $\varphi$ of the language of $\mathfrak{A}$, and for all $a_1, a_2, ..., a_n \in A$, $\mathfrak{A} \models \varphi(a_1, a_2, ..., a_n)$ iff $\mathfrak{B} \models \varphi(a_1, a_2, ..., a_n)$.} Thus, the resulting non-trivial reformulation of Self would, now, state that there is an elementary embedding $j:V\to M$, where $M$ is transitive and $j$ is not the identity, for then the pointwise image $V'$ of $V$ under $j$ is an isomorphic copy of $V$, different from $V$, which is an elementary substructure of $M$.

As a consequence, a way to express Self is given by elementary embeddings of the form: 

\[j: V \rightarrow M \]

\noindent
that are not the identity, and where $M$ is some transitive  subclass of $V$; so the image of $V$ under $j$ is an elementary  substructure of $M$, which can, finally, be conceptualised as a non-trivial `copy' of $V$. 

As is known, if such an embedding exists, then its critical point, that is, the least $\kappa$ such that $j(\kappa)\neq \kappa$ is a \textit{measurable} cardinal. The degree of closeness of $V$ and $M$ can be aptly calibrated by putting further closure constraints on $M$, e.g., by requiring that $M$ contains some large initial segment of $V$, or that $M$ is closed under long sequences. The relevant embeddings then give rise to different `copies' of $V$ that resemble $V$ more and more; in turn, this results in the postulation of the existence of ever stronger large cardinals, namely the critical points of the embeddings.  The upper bound for the template is given by Kunen's Theorem in \cite{kunen1971}, showing (in $\mathsf{ZFC}$) that $M$ cannot be $V$ itself. Moreover, \cite{suzuki1999} shows  that if the embedding $j$ is required to be definable, then this also holds in $\mathsf{ZF}$.

\medskip

Once made clear, through an analysis of Self, why one should use elementary embeddings to express Resemblance, one wonders whether this may also help us make the case that Resemblance (and Resemblance-inspired mathematical principles) are intrinsically justified (that is, whether Self is, itself, intrinsically justified). 

In connection with this goal, \cite{martin-steel1989} makes the following argument: if $\kappa$ is the critical point of an elementary embedding $j: V \rightarrow M$, with $M$ transitive, and, for some property $P$, $P(\kappa)$ holds in $V$, and $M$ resembles $V$ enough so that $P(\kappa)$ also holds in $M$, then since, necessarily, $\kappa < j(\kappa)$, it is true in $M$ that some $\alpha<j(\kappa)$ (namely $\kappa$), has the property $P$, and therefore, by the elementarity of $j$, there is also an $\alpha<\kappa$ which has the property $P$ in $V$. So, via $j$, the typical Reflection argument has gone through that, if a property $P$ holds in $V$ (of $\kappa$), there is some $\alpha<\kappa$ that already satisfies it.\footnote{\cite{martin-steel1989}, p. 73.} 

The argument showcases the connection between Reflection and (Self-moti-vated) Resemblance. \textit{Prima facie}, it seems to rest upon the assumption that the classes $j$ and $M$ are determinate entities, so, in the end, it might stumble upon the difficulties with proper classes we have already abundantly discussed. Nevertheless, in the presence of the Axiom of Choice, if an embedding $j:V\to M$ as above exists, then one can show, using an  ultrapower construction, that there is one with the same critical point $\kappa$ which is definable  using an ultrafilter over  $\kappa$ as a parameter.  

One could even be more audacious and state the existence of a direct link between (Self-motivated) Resemblance and Reflection through Unknowability, as suggested by \cite{hauser2006}, which countenances that elementary embeddings establish a connection with the Absolute.\footnote{\cite{hauser2006}, p. 536, cf. especially fn. 27.} But for this, clearly, one has to avail oneself of resources that are bound to overtake IC. 

All this suggests two tentative conclusions. On the one hand, the general template of elementary embeddings, hence, also Self, which is their underlying motivation, could, in principle, successfully meet Schindler's Challenge, insofar as elementary embeddings may be defined without involving arbitrary (impredicative) classes. This argumentative line will be further pursued through considering Structural Reflection (and its justifiability) in section \ref{SRP} and \ref{secfinal}. 

On the other, the mentioned connection between Reflection and elementary embeddings (Self), to the extent that the connection rests upon the `strong' mathematical form of the Reflection Principle discussed in section \ref{Refl}, does not, \textit{per se}, render Self (and elementary embeddings) justifiable in the light of IC, so, through Self, one still does not have at hand arguments establishing the intrinsic justifiability of Resemblance.  

\medskip

One final comment is in order. By especially looking at Self, one sees that Resemblance is a confluence of both Reflection and Uniformity ideas, and is, in a sense, a strengthening and an improvement on them. In particular, mathematical principles expressing Resemblance can be very strong, as they are able to capture very strong large-cardinal notions; so, these may also give it their best shot to provide a positive response to the Universality Issue. 

In what follows, we will review two mathematical expressions of Resemblance, and we will carefully explicate how these principles fail to solve the Intrinsicness Issue; a third one, Structural Reflection, will be discussed in a separate section, \ref{SRP}.

\subsubsection{Reinhardt's S4 and Welch's GRP}\label{RR}

\cite{reinhardt1974b} presents a class/set theory called $\mathsf{ZA}$ (after Zermelo and Ackermann). The language of $\mathsf{ZA}$ is the same as that of $\mathsf{ZF}$, but also contains individual constants for $V$ and $\On$; moreover, $\mathsf{ZA}$'s axioms are those of $\mathsf{ZF}$ (relativised to $V$) plus the additional axiom $V=V_{\On}$ and the Reflection Schema S2 discussed below. 

$\mathsf{ZA}$ is one further expression of \textit{set-theoretic potentialism} (see end of section \ref{IC}): in particular, Reinhardt takes $V$ to be a provisionally fixed universe, which can be extended to taller universes containing $V$ as a proper part, that he evocatively calls `projections of $V$'. To this end, in $\mathsf{ZA}$, one fixes the height of $V$ at $\On$ and can, then, define such objects (projections) as $V_{\On+1}$, $V_{\On+2}$, ...: all sets are in $V$, whereas `projections of $V$' contain `imaginings', sets which could conceivably be taken to exist (these are the \textit{classes} of $\mathsf{ZA}$).\footnote{\label{reinhardt}The bones of this conception are laid out in \cite{reinhardt1974b}, where Reinhardt introduces the distinction between `existing' and `imaginable' sets. For a brief, but exhaustive, discussion of the evolution of Reinhardt's ideas see \cite{kanamori2009}, pp. 312ff.} $\mathsf{ZA}$'s Reflection Principle S2 is the following axiom schema:
 
\[(\forall x,y \in V)[\theta^{V}(x,y) \leftrightarrow \theta(x,y)] \tag{S2} \] 

\noindent
where $\theta$ is any formula in the $\in$-language with free variables $x$ and $y$, and `$\theta^{V}$' is the relativisation of the quantifiers of  $\theta$  to $V$. Informally, S2 expresses the fact that a sentence about sets is true in $V$ if and only if it is true in a projection of $V$.

Now, Reinhardt proceeds to formulate, always in $\mathsf{ZA}$, the following stronger Reflection (in fact, Resemblance) \textit{schema}:

\[(\forall x,y \in V)(\forall P \subseteq V)[\theta^{\mathcal{P}(V)}(x,y, P) \leftrightarrow \theta(x,y, jP)] \tag{S4} \]

\noindent
where  $j$ is a unary function symbol, and $jP$ is the extension of $P$ in the projected universe beyond $V$.\footnote{Note that $j(x)=x$, if $x$ is a set, whereas in general $j(P) \neq P$, as $P$ is in $V_{\On+1}$, whereas $jP$ belongs to  $V_{j(\On)}+1$.} 

In essence, S4 further extends S2 by asserting that the properties of \textit{subclasses} of a provisionally fixed universe $V$\footnote{Expressed by sentences which may quantify over classes, and which may have sets and classes as parameters.} also hold of \textit{subclasses} of a projected universe beyond $V$,  whenever the class parameters $X$ are given their proper extensions $jX$, and vice versa. The quick reason for this is what Reinhardt calls the \textit{universality of set theory}, the principle whereby any entities beyond sets (if such things exist) should behave exactly like (should have the same properties of) sets.\footnote{\cite{reinhardt1974a}, p. 197.} We shall go back to this later in this section. 

In $\mathsf{ZFC}$, S4 could be aptly expressed in terms of the existence of an elementary embedding: 

\[j: V_{\alpha_{0}+1} \rightarrow V_{\alpha_{1}+1} \] 

\noindent
with $crit(j)=\alpha_0$; if $j$ exists, then $\alpha_0$ is a \textit{1-extendible} cardinal. Here, of course, $\On$ is interpreted as $\alpha_0$, $V$ as $V_{\alpha_0}$, $j\On$ as $\alpha_1$, and $jP$ as $j(P)$, for each $P\subseteq V_{\alpha_0}$.

\medskip

We now turn to examine \cite{horsten-welch2016}'s Global Reflection Principle (GRP). The underlying motivation for GRP is, in a sense, opposite to that for S4: if the latter is based on a \textit{potentialist} construal of $V$, the former presupposes a heavily \textit{actualist} interpretation of $V$, whereby not only $V$, but also its proper classes are fully determinate objects, a view the authors call Cantorian as opposed to the Zermelian (again, \textit{potentialist}) conception of $V$, since they take it to be an elaboration of Cantor's doctrine of the Absolute.\footnote{Cf. \cite{horsten-welch2016}, especially section I. \cite{barton2015} provides an account of GRP based, on the contrary, on a Zermelian conception of $V$, and on `richness' considerations along the lines of those in section \ref{VP}. However, Barton's interpretation comes at the expense of no longer viewing GRP as a `stable' set-theoretic truth (cf. \cite{barton2015}, pp. 356-7).}

Sets in $V$ and proper classes, are, then, taken to be reflected downwards by, respectively, the sets of a rank initial segment $V_{\kappa}$ of $V$ and by a collection of subsets of $V_{\kappa}$, through positing the existence of an elementary embedding:

\[j: (V_{\kappa}, \in, V_{\kappa +1}) \rightarrow (V, \in, \mathcal{C}) \tag{GRP} \]

\noindent
where   $\mathcal{C}$ is the collection of all subclasses of $V$; so, GRP is precisely the statement that such a $j$ exists. Although GRP, thus formulated, is a third-order statement, as $\mathcal{C}$ is a third-order predicate over $V$, \cite{welch2017} shows that it can actually be formulated so that $j$ is a  second-order object, and truth in $(V, \in, \mathcal{C})$ can also be formulated as a second-order relation by adding a second-order satisfaction predicate to the language. The \emph{elementarity} of the embedding $j$ is understood in the usual sense in the language of first-order set theory enriched with a collection of predicate variables, $X_0$, $X_1$, $X_2$ (varying over $\mathcal{C}$), but quantification over these is barred: this allows the authors to claim that GRP does not commit one to the existence of a \textit{completed} realm of classes (and to higher-order quantification). More on this in a moment.

GRP implies that $\kappa$ is a measurable  Woodin cardinal, and, moreover, that there are \textit{unboundedly} many measurable Woodin cardinals beyond $\kappa$.\footnote{\cite{welch2017}, p. 99.  A cardinal $\kappa$ is Woodin if and only if for any $f: \kappa \rightarrow \kappa$, there exists an $\alpha<\kappa$ which is closed under $f$, and an elementary embedding $j: V \rightarrow M$, whose critical point is $\alpha$, such that $V_{j(f)(\alpha)} \subseteq M$.} As a consequence, GRP also implies other fundamental set-theoretic statements, including PD, AD$^{L(\mathbb{R})}$, and that the theory of $L(\mathbb{R})$ is absolute among set-generic extensions of $V$.\footnote{\label{welch}\cite{welch2017}, p. 99.} In terms of consistency strength, the consistency of the existence of a $1$-extendible cardinal yields the consistency of $\mathsf{NBG}$ plus GRP.\footnote{\cite{welch2017}, Lemma 8.} However, $\mathsf{NBG}$ plus GRP does not directly yield the existence of a $1$-extendible cardinal, although raising the elementarity of the embedding $j$ to $\Sigma^1_1$-elementarity does imply it.\footnote{\cite{horsten-welch2016}, Section IV.} \cite{welch2017} also considers further strengthenings of GRP, such as GRP$^{1}_{\infty}$, arising from lifting the level of elementarity of $j$ in GRP to $\Sigma^{1}_{\infty}$-formulas, and GRP$^{+}$, arising from enhancing the language of second-order set theory with a satisfaction predicate. \cite{roberts2017}'s theory $\mathsf{ZFC}_{2}+R_{S}$ (introduced in section \ref{Refl}) also follows suit (although its reflection principle $R_{S}$ is, in the face of it, different from GRP), and is, as a consequence, able to yield a \textit{proper class} of 1-extendible cardinals.\footnote{\cite{roberts2017}, p. 655. For the relationship between \cite{roberts2017}'s principles and GRP, also see \cite{welch2019}, p. 95.} 

\medskip

\begin{remark}\label{remark1}

\normalfont S4 and GRP are not the only existing mathematical formulations of Resemblance. We have briefly mentioned \cite{roberts2017}'s principles, which, however, might be taken to represent generalisations (strengthenings) of GRP. \cite{marshall1989}'s theory $B_2$ aims to be a generalisation of \cite{bernays1976}'s theory $B_1$ featuring ${\rm RP}\Pi^1_1$. To this end, Marshall takes into account 2-classes, that is, classes of classes, and defines a Reflection Principle (A3) that applies to 2-classes in the same way as ${\rm RP}\Pi^1_1$ applied to (1-)classes. $B_2$ is very strong, insofar as it proves the existence of a proper class of \textit{1-extendible} cardinals. In a similar vein as \cite{reinhardt1974a}'s account of $\Omega$-classes, then Marshall proceeds to define a sequence of stronger theories with $n>2$-classes, and extensions of such theories are able to capture even stronger large cardinals. For the potential connections between this approach and Uniformity, see below Remark \ref{reinhardt-uniformity}.\footnote{For the connections between \cite{roberts2017}'s and \cite{marshall1989}'s theories, also see \cite{mccallum2021}, pp. 207ff.}    
    
\end{remark}

\medskip

\noindent
Let's take stock. We have examined two mathematical expressions of Resemblance, S4 and GRP. We should now consider how they fare with the Intrinsicness Issue.

Reinhardt has provided several conceptual underpinnings for S4. One is the mentioned idea of the `universality of set theory'; a second one, we will not take into account, even adumbrates the use of theological principles, in particular the Pseudo-Areopagite's doctrine of analogy.\footnote{Cf. \cite{reinhardt1974a}, Remark 5.6, p. 198.} A third one is based on Reinhardt's own potentialist 
conception of collections, some of which, the \textit{sets}, are taken to be actual, others (\textit{classes}) are just `imaginable' objects.\footnote{Cf. fn. \ref{reinhardt}. A full modal treatment of sets and intensional properties is to be found in \cite{reinhardt1980}.}  

Now, a more determinately actualist intepretation of S4 may also be provided. Based on the G\"odelian considerations in \cite{wang1996}, one could see Reinhardt's S4 as being naturally motivated by a strengthening of Unknowability reformulated as follows: 

\vspace{11pt}

\noindent
\textbf{Unknowability}. The Absolute, \textit{as well as its properties} (in an idealised sense), are unknowable. 

\vspace{11pt}

In particular, in \cite{wang1996}, G\"odel suggests that S4 (in fact, its precursor, the axiom S3.3 of \cite{reinhardt1974b}) would be a natural strengthening of Ackermann's Axiom $\gamma$, a Reflection Principle itself.\footnote{\label{ackermann}Cf. \cite{wang1996}, p. 285; \cite{wang1977}, p. 325. \cite{ackermann1956}'s Axiom $\gamma$ is: $(\forall x)(\mathfrak{A}(x)\rightarrow M(x))\rightarrow(\exists y)(z \in y \leftrightarrow \mathfrak{A}(z))$ [If there is a condition such that only sets satisfy it, then there exists the set of all sets which satisfy it]. `$M(x)$' means `$x$ is a set', $\mathfrak{A}(x)$ is a property that may also contain set parameters, but, crucially, may not contain the predicate `$M(x)$'. It should be noted that, in Ackermann's view, Axiom $\gamma$ was intuitively motivated by the idea that what sets \textit{there are} cannot be etsablished \textit{a priori}. But in Ackermann's (as well as Reinhardt's) theory, this resulted in taking classes such as $V$ to be \textit{indeterminate} in a \textit{potentialist} sense, a view which does not square well with G\"odel's \textit{actualist} interpretation of the Absolute and of Unknowability we present in this section.} 
Now, since, according to G\"odel, Ackermann's Axiom $\gamma$ would be justified by (would even be equivalent to) Unknowability, Reinhardt's S4 could be seen as also being justified by Unknowability (in the form stated above), since, now, S4 would reflect (be motivated by) the fact that it is not only $V$ to be unknowable, but also its properties, taken to correspond, logically, to \textit{proper classes}.

\medskip

Let us consider these strategies in connection with the Intrinsicness Issue. The potentialist one would seem to be compatible with IC, insofar as it does not commit itself to an \textit{actualised} $V$. Moreover, the `universality-of-set-theory argument' guarantees that any entity of set theory, beyond sets themselves, is, to all intents and purposes, reducible to \textit{sets}, that is, to the entities formed according to IC. However, as we have seen, $\mathsf{ZA}$ commits itself to classes (second-order objects) a lot more than it would seem at first glance: $\mathsf{ZA}$, in particular, \textit{impredicatively} quantifies over them, hence needs to take classes to constitute a \textit{fully determinate} realm of objects; moreover, for the elementary embedding $j$ in S4 one needs to `move' classes to other classes in a way which cannot be circumscribed by any definition.\footnote{For the latter criticism, see \cite{koellner2009b}, p. 217.}    

Also, if one paraphrases S4 as asserting the \textit{indescribability} of the Cantorian Absolute, that is, as just more strongly and determinately expressing Unknowability, in accordance with an actualist conception of $V$, then one, again, needs to refer to the Absolute, and, hence, use resources which overtake those of IC. Either way, one runs into damning difficulties. 

\begin{remark}\label{reinhardt-uniformity}

\normalfont

S4 can be strengthened by taking into account 3rd-order, 4th-order, $n$-th order classes, etc. \cite{reinhardt1974b} globally calls these $\Omega$-classes. The same construction is hinted at in \cite{marshall1989} (cf. Remark \ref{remark1}). Now, Reinhardt posits that, for all $\lambda \in \On$, if $\lambda$ is an $\Omega$-class, $\Omega$ is $\lambda$-extendible (this is \cite{reinhardt1974a}'s Axiom 6.3). Thus, mathematically, through Axiom 6.3, one gets \textit{extendible} cardinals. But now how can one justify talk of `$\Omega$-classes' and the extension of Reflection to all of these? One strategy consists in further modifying Unknowability so as to include mentioning of \textit{higher-order} properties. Another strategy consists in viewing this move as being justified by G\"odel's Uniformity as discussed in section \ref{uniformity}: `properties of sets re-appear \textit{over and over} in $V$', where `over and over in $V$' is taken care of by the extension of the relevant elementary embeddings $j$'s that define the $\lambda$-extendible cardinals through the hierarchy of $\Omega$-classes. The resulting Reflection Principles (analogous to S4) would, now, quantify over an infinite  hierarchy of new objects, $\Omega$-classes, none of which, as is clear, is justified by IC. 
    
\end{remark}

The relationship between GRP and the Intrinsicness Issue can be more easily assessed. \cite{horsten-welch2016} takes classes to be construable as \textit{parts} of $V$, that is, as entities distinct from sets, but existing \textit{alongside} $V$ \textit{mereologically}, a conception indebted, as the authors make it clear, to \cite{lewis1991}. Now, regardless of whether this interpretation of GRP is tenable, the authors explicitly trace their conception of classes, as we have seen, to Cantor's Absolute, thus hinting at the existence of a strong conceptual link between GRP, the Absolute, and the collection of \textit{all} proper classes of $V$ (the hyper-class $\mathcal{C}$ in GRP). As a consequence, Cantor's Absolute should be taken to be indispensable to the correct understanding and even formulation of GRP.

As we have seen in the previous paragraphs, the authors lay considerable emphasis on the fact that the elementarity of $j$ in the statement of GRP is understood in the first-order sense (i.e., $\Sigma^0_\infty$-elementarity in the language expanded with new  variables that range over $\mathcal{C}$): it is only in the stronger principles, such as  ${\rm GRP}_{\Sigma^{1}_{1}}$, which requires $\Sigma^1_1$-elementarity,  that quantification over arbitrary classes is allowed, thus, formally, overtaking the boundaries of IC. For that matter, the authors also make it clear that, if one wishes, one could formulate GRP in $\mathsf{NBG}$, but that this, strictly speaking, is not necessary, as they think of the structure $(V, \in, \mathcal{C})$ as enjoying a `pre-formalised' state.\footnote{Cf. \cite{welch2019}, p. 94.} 

In any case, given the strategic importance of the Absolute (and of its `parts') as an intuitive rationale behind GRP, it seems fair to say that the principle does, indeed, need more conceptual resources than those afforded by IC in order to be properly motivated.\footnote{Incidentally, \cite{roberts2017}, pp. 657-8, also shares this view concerning the status of its own $R_{S}$.} 

\medskip

In sum, both S4 and GRP do not seem to be able to solve the Intrinsicness Issue, or, for that matter, meet Schindler's Challenge. In particular, they make use of either arbitrary proper classes, or analogous concepts, such as that of `parts', not licensed by IC. Incidentally, they also look unnatural, as they require considerable tweaking and reinterpretation of the natural notions of `set', `class', and, even, 'elementary embedding'. Hence, they are somewhat contrived in their formulations, and therefore hard to use.  

As far as the Universality Issue is concerned, both S4 and GRP are very strong Resemblance principles, as they may even be able to yield \textit{extendible} cardinals (as far as GRP is concerned, it is only its strengthening ${\rm GRP}_{\Sigma^{1}_{1}}$ which is able to yield them). So, in principle, they may be able to positively respond to the issue. However, even if they really do that, they do that \textit{unsystematically}, so to speak, without proposing a transparent and workable methodology.

\section{Structural Reflection}\label{SRP}

Structural Reflection will be our main focus in the second part of the paper. We first formulate it, and then discuss its features, strengthenings and consequences. Its basic form is the following schema \cite{bagaria2023}:

\begin{definition}[Structural Reflection (SR)]

For every definable, in the first-order language of set theory, possibly with parameters, class $\mathcal{C}$ of relational structures of the same type there exists an ordinal $\alpha$ that \textit{reflects} $\mathcal{C}$, i.e., for every $A$ in $\mathcal{C}$ there exists $B$ in $\mathcal{C} \cap V_{\alpha}$, and an elementary embedding from $B$ into $A$. 

\end{definition}

SR arises from the confluence of several notions and mathematical principles, some of which we have already examined; most prominently  the following three:

\begin{enumerate}

    \item G\"odel's notion of `[reflecting an] internal structural property of the membership relation in $V$'

    \item Vop\u{e}nka's Principle (see section \ref{VP})

    \item Resemblance (see section \ref{self})
    
\end{enumerate}

Some observations are in order. As far as (1.) is concerned, \cite{bagaria2023} explicitly construes G\"odel's notion of `internal structural property of the membership relation in $V$', discussed in \cite{wang1996},\footnote{Cf. \cite{wang1996}, pp. 283-4. Wang says that G\"odel took structural properties to be those used in Ackermann's Axiom $\gamma$ (for which see fn. \ref{ackermann}), i.e., all those properties (expressible in the language of $\mathsf{A}$) which do not explicitly mention $V$ (this being an instance, in G\"odel's view, of Unknowability).} as being expressed by a \textit{definable class} of relational structures (of the same type) $\langle A, \in, \langle R_{i} \rangle_{i \in I}\rangle$, where $A$ is a nonempty set, and $R_{i}$, $i\in I$, is a family of relations on $A$. 

As for (2.), a quick comparison  shows that SR implies Vop\u{e}nka's Principle, formulated as a schema; as has already been said, the latter states that, in any definable proper class of structures of the same type there are \emph{two} members of the class, $A$ and $B$, such that $B$ is elementarily embeddable into $A$, whereas SR states that, for \textit{every} $A$ in the class, \textit{there is} a $B$ also in the class, but of rank lower than some fixed rank $\alpha$, such that $B$ is elementarily embeddable into $A$. In fact, SR and VP, both taken as schemata, are equivalent (Theorem \ref{thmSRequivVP} below). Hence, the Richness considerations that \cite{maddy1988b} saw as motivating Vop\u{e}nka's Principle are also applicable to SR. This is of some relevance to our arguments in section \ref{just}.

With regards to (3.), SR is, like the other Resemblance principles, also crucially hinged on the use of elementary embeddings. However, in the case of SR the elementary embeddings are always \textit{between sets}, more precisely, between structures of the same type, which in most cases may be assumed to be of the form $<V_\beta ,\in, A>$, with $A$ either an element or a subset of $V_\beta$. As the goal of SR is to state that \textit{every} structural property of  $V$, as conveyed by a definable class of relational structures, is reflected in some fixed rank-initial segment $V_\alpha$ of $V$, the existence of an elementary embedding from some small structure in the class, namely one that belongs to $V_\alpha$, into any given  structure in the class, which possibly has much larger rank or size, yields a precise mathematical formulation of the intuitive notion of resemblance between the two structures. Let us also emphasise that although, as we have seen, other Resemblance principles also take \textit{classes} to express `properties of $V$', SR differentiates itself from those principles in that it considers only \emph{definable} classes, i.e., classes defined  by \textit{formulas} of the language of set theory, possibly with sets as  parameters.\footnote{Cf. \cite{bagaria2023}, pp. 29ff.}  

\medskip

As said, SR is just the general form  of a   family of Structural Reflection Principles (SRPs). This is because, to begin with, SR cannot be formulated in the first-order language of set theory as a single statement, as to \textit{each} definable class $\mathcal{C}$ there corresponds some first-order formula in the language of set theory, and those formulas have unbounded complexity. Thus, SR must be formulated as a schema, namely as an infinite list of statements that are obtained by restricting the  formulas involved to those that have bounded complexity, namely to  $\Sigma_{n}$- or $\Pi_{n}$-formulas in the Lévy hierarchy, for each natural number $n$. As it turns out, SR restricted to $\Sigma_{n+1}$-definable classes is equivalent to SR restricted to $\Pi_n$-definable classes, so SR may be formulated as  the schema asserting $\Pi_n$-SR, for each $n$:

\begin{definition}[$\Pi_n$-SR]

For every $\Pi_n$-definable class $\mathcal{C}$ of structures of the same type there is an ordinal $\alpha$ that \textit{reflects} $\mathcal{C}$.

\end{definition}

The principles $\Pi_n$-SR are just the simplest members of a complex and ramified network of Structural Reflection Principles obtained through introducing further constraints and conditions, or by further generalising them. The articles \cite{bagaria2023}, \cite{bagaria-vaananen2016}, \cite{bagaria-lucke2022}, \cite{bagaria-wilson2023}, and \cite{bagaria-lucke2023}, describe many other types of SR, including: SR$^{-}$, Product Structural Reflection (PSR), Structural Reflection restricted to \textit{inner models}, Generic Structural Reflection (GSR), Weak Structural reflection (WSR), Exact Structural Reflection (ESR), etc. Work on producing variants of basic SR has been propelled by, among other things, the attempt to make SR as \textit{universal} as possible, that is, to enable it to capture as many large-cardinal notions as possible across the whole spectrum of consistency strengths of large cardinals. 

\medskip

In its basic form, SR  is very strong, as shown by the following: 

\begin{theorem}[\cite{bagaria2023}, p. 32] \
\label{thmSRequivVP}
\begin{enumerate}
    \item $\Pi_1$-SR  is equivalent to the existence of a \textit{supercompact} cardinal.

    \item $\Pi_2$-SR  is equivalent to the existence of an \textit{extendible} cardinal. 

    \item SR is equivalent to VP (taken as a \textit{schema}).\footnote{Originally proved in \cite{bagaria2012}, pp. 223ff. With the \emph{boldface} versions of these, i.e., by allowing parameters in  the formulas defining the classes of structures, one obtains proper classes of the corresponding large-cardinal notions; for instance, $\boldsymbol{\Pi_{1}}$-SR is equivalent to the existence of a \textit{proper class} of supercompact cardinals. However, the two versions of VP, the lightface and the boldface, are equivalent.}

\end{enumerate}

\end{theorem}
  
As a general taxonomic fact, each type of SRP will specialise, as it were, in a particular region of the consistency-strength hierarchy of large cardinals.\footnote{Cf. the summarising tables at the end of \cite{bagaria2023}, pp. 65-7.} For instance:

\begin{itemize}
    \item $\Pi_n$-SR, for $n>0$, characterises large cardinals from supercompact to extendible, up to the level of Vop\u{e}nka's Principle;

    \item SR$^{-}$ yields characterisations of  small large cardinals (inaccessible, Mahlo, weakly compact);

    \item PSR (Product SR) yields characterisations of large cardinals from strong to Woodin cardinals;

    \item GSR (Generic SR) and Strong GSR yield characterisations of large cardinals from almost remarkable to virtually extendible cardinals;

    \item WSR (Weak SR) characterises large cardinals between strongly unfoldable and subtle;

    \item ESR (Exact SR) characterises large cardinals beyond Vop\u{e}nka's Principle, up to the level of I1 embeddings;

\end{itemize}

This complex, and multi-faceted, ramification of principles, and the equivalence results with large-cardinal notions obtained, has recently led the first author to formulate the following conjecture:

\begin{quote}

Each of these results should be regarded as a small step towards the ultimate objective of showing that all large cardinals are in fact different manifestations of a single general reflection principle. \\(\cite{bagaria2023}, p. 30)

\end{quote}

\noindent
that is:

\begin{conjecture}\label{conj}

Every known large cardinal notion is equivalent to some (form of) SRP.
    
\end{conjecture}

\noindent
In the light of Conjecture \ref{conj}, one could formulate one further natural  conjecture, which attempts to define what is (the least cardinal satisfying) a large cardinal notion, which might count as a fully general definition of \textit{necessary} and \textit{sufficient} conditions for some cardinal $\kappa$ to be a \textit{large} cardinal:

\begin{conjecture}\label{conj2}

A cardinal $\kappa$ is the least cardinal satisfying some large-cardinal notion iff $\kappa$ is the least cardinal satisfying some Structural Reflection Principle that implies  (in some inner model) $\kappa$ is (weakly) inaccessible.
    
\end{conjecture}

We close this section with two further considerations about how SRPs further differentiate themselves from the other Resemblance principles. 

First, we wish to stress again that classes, in SRPs, are always \textit{definable} classes (with parameters), whereas, as we have seen, the other principles, very often, need to quantify over (all) classes. Secondly, contrary to the other principles, SRPs provide one with a systematic approach  to yield large-cardinal notions. Both aspects of SRPs will prove fundamental for our discussion of the justifiability of large cardinals through SR in section \ref{just}.

\section{Large Cardinals Under Structural Reflection}\label{secfinal}

\subsection{Justifying Large Cardinals Through SRPs}\label{just}

We are now ready to consider how SRPs fare with the Issues we have discussed throughout the paper and, in particular, with the intrinsic justifiability of large cardinals. 

We have seen that SR is a Resemblance principle. However, since it also provides a mathematical expression for G\"odel's notion of reflection of a structural property,   SR is also, conceivably, motivated by Reflection. Let us look into this in more detail. 

It seems fair to say that G\"odel thought Reflection to be equivalent to Unknowability (in particular, to Unknowability as stated in section \ref{RR}). This means that, since, on G\"odel's view, $V$ is indefinable (unknowable), every structural property of the membership relation  should be taken to be a property of that relation in some initial segment of $V$. So, one could say that Unknowability, in particular, motivates SR, insofar as the latter prescribes that \textit{definable} classes of relational structures, the embodiment of G\"odel's notion of structural property,\footnote{See \cite{bagaria2023}, section 2.} are reflected in initial segments of $V$. However, this motivation will be of little to no effect on the intrinsic justifiability of SR, since, as we know from section \ref{Refl}, with its essential reference to \textit{completed} totalities, Reflection uses resources which go beyond IC.  However, as we shall see in a moment, Reflection, in any of its  forms that make reference to classes as completed totalities, does not actually need to be woven into the justification of SRPs. 

Let us, then, proceed to consider Resemblance. We have examined two  mathematical expressions of this principle: Reinhardt's S4 and Welch's GRP. Both assume the existence of a domain of classes which need not be restricted to those that are definable. 
In contrast, the classes that SRPs use are only the \textit{definable} ones. This, among other things, aligns SRPs with the practising set-theorist's presupposition that classes, in $\mathsf{ZFC}$, are just given by formulas of the language of set theory, with   parameters. Now, formulas are not part of the presupposed realm of objects (ontology) of $\mathsf{ZFC}$, or, if they are so (e.g., if one thinks of them  as natural numbers, via some coding), they are still sets, so SRPs do not presuppose the existence of other objects than sets. 
However, one could object that SRPs are hinged on elementary embeddings. As we have seen in section \ref{self}, in the usual formulation of large cardinals such as measurable, strong, supercompact, etc., elementary embeddings are of the form $j: V \rightarrow M$, where $M$ is a transitive class. So they are \textit{class-embeddings} (of the proper class $V$ into a subclass $M$ of $V$ itself).
But in the formulation of SRPs, the embeddings in use are only set-embeddings: in the statement ``for some $A$ and $B$ in $\mathcal{C}$, there is an elementary embedding $j: A \rightarrow B$'', $A$ and $B$ are \textit{sets}, so $j$ is just another set. 

All this stands SRPs in good stead with regard to Schindler's Challenge: SRPs are reflection principles yielding very strong large cardinals, yet do not assume the existence of impredicative classes. But does this, crucially, count as a solution to the Intrinsicness Issue? Our answer is: yes and no. 

In a sense, yes, since Schindler's Challenge precisely challenges the fact that one may formulate principles equivalent to large cardinals without, thereby, also positing the existence of classes that  are not licensed by IC; now, work on SRPs has, precisely, shown that this need not be the case, that is, that the formulation of SRPs does not require resources that overtake the boundaries of IC. 

In an another sense, no, since the resources used by SRPs crucially include \textit{elementary embeddings} (although just set-embeddings) for whose justification one may either invoke Richness-based considerations like those used by Maddy to support VP (section \ref{VP}), or Resemblance-based considerations (section \ref{self}), particularly, those relating to the principle we have called Self; in either case, one runs into already amply discussed difficulties.  

For another issue, some SRPs we have reviewed in section \ref{SRP}, e.g., $\Gamma$-SR, Product SR, Generic SR, etc., seem to be even less amenable to IC's resources, insofar as they contain references to purely model-theoretic constructions, such as  \textit{product of structures}, \textit{generic filters}, and others. Although the latter are, in the end, just other set-theoretic entities, they are not, \textit{per se}, straightforwardly amenable to the rationale behind IC. However, as a potential mitigating circumstance, one should note that `simple' Structural Reflection may be immune to this objection, as we know (section \ref{SRP}) that that form of SR is equivalent to very strong large-cardinal notions (even to VP). If SR is, then, able to justify VP intrinsically, then one could argue that all  weaker large cardinals are, \textit{a fortiori}, intrinsically justified by SR. The only SRPs which would, thus, be prone to the objection are those yielding large cardinals stronger  than VP; hence, by this argument's lights, it would only be ESR (again, see section \ref{SRP}), and the corresponding large-cardinal notions, that would be lacking intrinsic evidence. 

The Universality Issue is an entirely different story, since, as we have seen, the results attained so far, some of which have been reviewed in section \ref{SRP}, even point to a validation of Conjecture \ref{conj}. Moreover, by coming with a carefully thought out methodology, which is liable to be further refined and strengthened, to accommodate all kinds of large-cardinal notions, contrary to S4, GRP and other Resemblance principles, SR is also able to \textit{explain} why it can provide, in a non-haphazard manner, a positive response to the Universality Issue.  

\medskip

It seems to us that we have now arrived at a crucial crossroads: if Conjecture \ref{conj} is validated, and, in addition, SRPs are intrinsically justified, even Koellner's Challenge (cf. section \ref{Refl}) is met by SRPs; if SRPs are not intrinsically justified, given how well they fare with the Universality Issue, they could, regardless, be taken to be the right extension of $\mathsf{ZFC}$. In the latter case, however, SRPs' striking ability to provide a positive response to the Universality Issue, may, in the end, in a way or the other, also reflect on their intrinsic justifiability. We shall say more on this in section \ref{equilibrium}. 

\medskip

There is one last strand of results which is worth mentioning at this stage, in light of which SR seems to be justified, results which belong, this time, to the realm of pure logic. 

Recall that the L\"owenheim-Skolem Theorem  (LST) for first-order logic asserts that every structure, in a countable language, has a countable  elementary substructure. More generally, every structure, in a language of size $<\!\kappa$, where $\kappa$ is a regular uncountable cardinal, has an elementary substructure of size $<\!\kappa$. Thus, SR is just a strong form of LST, for it asserts that for any definable property $P$ of structures in the same language, there is a regular uncountable cardinal $\kappa$, greater than the size of the language, such that every structure with the property $P$ has an elementary substructure of size $<\!\kappa$ which is isomorphic to a structure that also has the property $P$. Thus, SR is just like LST with the additional requirement that the small elementary substructure (is isomorphic to a structure that) has a prescribed property $P$. We will briefly go back to the connection between SR and LST when discussing the justification of SR in section \ref{equilibrium}.

\medskip

\begin{remark}\label{LST}

\normalfont

It has long been known that strong forms of LST yield large cardinals. In particular, the following theorem of Magidor shows that the least supercompact cardinal is precisely the \textit{least} cardinal $\kappa$ that witnesses LST for second-order logic. Namely:

\begin{theorem}[\cite{magidor1971}]\label{mag}

A cardinal $\kappa$ is the least supercompact cardinal if and only if it is the least cardinal such that for every second-order formula $\varphi$, and every structure $\mathcal{A}$ in the language of $\varphi$, if $\mathcal{A} \models \varphi$, then there is a substructure $\mathcal{B}$ of $\mathcal{A}$ of cardinality less than $\kappa$ such that also $\mathcal{B} \models \varphi$. 
    
\end{theorem}

Now, if we compare the characterisation of the first supercompact cardinal given by Theorem \ref{mag} with the characterisation given by SR, namely as the least cardinal that reflects the class of $V_\alpha$, or equivalently, the least cardinal that reflects every $\Pi_1$-definable class of structures (Theorem \ref{thmSRequivVP} above), we see that the difference is the following: second-order reflection of formulas is replaced in the SR characterisation by full first-order reflection, but the small substructure $\mathcal{B}$ of $\mathcal{A}$ is now required to be isomorphic to some structure that has some predetermined $\Pi_1$-definable property that $\mathcal{A}$ has (equivalently, the property of being a $V_{\alpha}$). Thus, if the $\Pi_1$ property is to be a $V_\alpha$, for some $\alpha$, and $\mathcal{A}$ satisfies this property, and so it is  a $V_\alpha$, then $\mathcal{B}$ must be also (isomorphic to) some $V_\beta$.

\end{remark}

\subsection{Classes \textit{vs.} IC}\label{interlude}

One could attempt to respond to the Intrinsicness Issue (and to Schindler's Challenge) in a slightly different way, that is, by arguing that classes are part of (are implied by) the concept of set. One could proceed as follows. First, one could question the fact that the concept of set is IC. Other 
(arguably, less used in practice) conceptions/explications of sets have been formulated over the years, such as the \textit{limitation of size} conception, the \textit{logical} conception, \textit{plural quantification}, varieties of \textit{potentialism}, and others, which might be preferable to IC and, in addition, might justify the existence of classes.\footnote{For an overview of `conceptions' (as opposed to `concepts') of set, cf. \cite{incurvati2020}.}.

There are two main problems with all these conceptions: first, even if they were able to persuade us of the conceptual tenability of classes, they wouldn't, \textit{per se}, make their \textit{use} in set theory more palatable to the working set-theorist; second, they might be prone to overgeneration, that is, they might be liable to license also the generation of further entities beyond classes themselves.\footnote{As \cite{fraenkel1973} once observed:

\begin{quote}
    This process of adding bigger classes and hyper-classes has to stop somewhere; and we have to decide where to do so. $\mathsf{QM}$ [$\mathsf{MK}$, our note] is a good place to stop at for reasons of convenience and neatness, yet, apart from these considerations, this choice is as arbitrary as any other. (p. 145)

\end{quote}.}

We will just touch on the first problem. A way to address it would be to have at hand an account of what classes are \textit{for} in set theory, but, since, to begin with, there is no widespread consensus about which account of classes is correct, the practical indispensability of classes remains very suspect.\footnote{Cf. \cite{maddy1983}, and \cite{fujimoto2019}, which discusses four main accounts of classes. A useful primer on class theories is in \cite{potter2004}, Appendices B and C.} As a consequence, also the extent of the justifiability of large cardinals through a full-blown class-theoretic approach greatly varies. E.g., the Zermelian\footnote{Cf. fn. \ref{zermelo}.} may not go so far as justify large cardinals stronger than strong inaccessibles; Ackermann's concept of class, through its extension in Reinhardt's $\mathsf{ZA}$, might be able to justify large cardinals up to the level of \textit{extendible} cardinals; on the other hand, \cite{horsten-welch2016}'s rich realm of \textit{mereological} classes is equivalent to, and helps justify, strong large-cardinal notions, such as as measurable and Woodin, but, perhaps, none stronger than these, and so on. 

\medskip

Whatever the correct interpretation of classes, their use in the theory of large cardinals might, in any case, be, as suggested by \cite{uzquiano2003}, `heuristically indispensable'.\footnote{\cite{uzquiano2003}, p. 70.} But then, for this, i.e., for the purpose of providing a useful heuristic, the ontologically austere approach we have adopted and advocated in section \ref{just}, \textit{definabilism}, seems to be wholly adequate.
    
Recall that \textit{definabilism} is the idea that classes are given by  formulas in the meta-theory, that is, that they do not exist independently of definitions. Now, \cite{fujimoto2019} has argued that \textit{definabilism} is not adequate to interpret many instances of set-theoretic discourse, where classes are invoked (and, in Fujimoto's view, are even \textit{formally} indispensable).

In particular, Fujimoto argues that definabilism will not meet two essential desiderata of a correct interpretation of classes: (i) non-trivialisation of mathematical results; (ii) use of an appropriate mathematical framework for the class-theoretic discourse itself. Now, the significance of both criteria is debatable. In particular, it would seem that both are too \textit{ad hoc} to be suitable to the issue under consideration, as they seem to express concerns to which the definabilist will, presumably, be insensitive.\footnote{This seems to transpire from the examples made by Fujimoto. With reference to (i), Fujimoto considers \cite{suzuki1999}'s  proof of Kunen's inconsistency (in $\mathsf{ZF}$) for \emph{definable} elementary embeddings $j:V\to V$, and argues that Suzuki's result is not regarded by set-theorists as being `equivalent' to Kunen's proof, which uses AC, and which holds for \emph{all} such embeddings. Note, however, that Kunen's inconsistency can just be viewed as showing that `there is no non-trivial elementary embedding $j:V_{\lambda +2}\to V_{\lambda +2}$', a statement which, by no means, involves the use of proper classes. As regards (ii), Fujimoto makes the example of S. Friedman's Inner Model Hypothesis (IMH) and $V$-logic, introduced in \cite{afht2015}, and then further discussed in \cite{antos-barton-friedman2021}. The author correctly points out that classes may be needed to code `extensions of $V$', in $V$-logic, but fails to mention that the Completeness Theorem of $V$-logic, whereby one has `real' extensions of $V$, only goes through if $V$ is countable, so the `countable-$V$-logic' approach to IMH might even turn out to be more mathematically appropriate than one based on taking $V$ to be class-sized. As a consequence, Fujimoto's example seems to be rather inconclusive.} 

Moreover, the fact that, as we have seen in section \ref{SRP}, SRPs may be able to yield \textit{all} large cardinals makes the case for definabilism even more solid and, in any case, any working set-theorist will be content with practising that form of virtuous Ockham's razor which consists in denying reality to classes despite acknowledging their usefulness (if any): the details, then, of exactly how much class theory is really needed to investigate large cardinals will just be seen as a matter of formal convenience by those set-theorists. In practice, one needs at most to expand the language of set theory by adding predicate symbols for proper-class elementary embeddings, but even in this case there are ways to avoid taking this step.

Finally, whenever pressed to construe set-theoretic discourse involving clas-ses in a more faithful way, the definabilist can always resort to a reduction of classes to sets through what \cite{fujimoto2019} himself calls `hermeneutic reductionism', i.e., the idea that statements about classes should be taken to refer to set-sized structures $\langle V_{\kappa}, \in, V_{\kappa+1} \rangle$, with $\kappa$ a strongly inaccessible cardinal.\footnote{\cite{fujimoto2019}, p. 207ff.} Although also the latter account is not unproblematic, it will invariably provide working set-theorists with the required mathematical understanding of classes needed in most set-theoretic contexts.
    
\medskip

For another, more extreme, approach, one could, indeed, take \textit{classes} to be the real subject of set theory, and \textit{sets} just a type of \textit{classes}. This is, for instance, the conception underlying Ackermann's set theory $\mathsf{A}$, and, to some extent, also Reinhardt's $\mathsf{ZA}$. More recently, \cite{muller2001} has also taken this route and formulated a \textit{class-theoretic reductionist} account of the foundations of mathematics that he describes as follows:

\begin{quote}
    Every mathematical object is a class that demonstrably exists; and every variable used in some branch of mathematics is a variable over a fixed range, which is a class that demonstrably exists. \\ (\cite{muller2001}, p. 554)
\end{quote}

Conceptions such as \cite{muller2001}'s revolve around the idea that the concept of \textit{class} (even of \textit{proper} class) is prior to (the concept of) \textit{set}. We have abundantly commented on the troubles with Reinhardt's (thus, indirectly also Ackermann's) kindred conception. 

With reference to our purposes in this subsection, the main problem with accepting this kind of class-theoretic reductionism is that it has no clear advantages over set-theoretic reductionism as far as large cardinals are concerned: in particular, it is not clear that M\"uller's class theory $\mathsf{ARC}$ is suitable to motivate principles sufficiently strong to produce all large-cardinal notions.\footnote{Especially given the equiconsistency of $\mathsf{ARC}$ with just the theory $\mathsf{ZFC}$+`there exists a strongly inaccessible cardinal'. Cf. \cite{muller2001}, pp. 563-4.}  

\medskip

So, in conclusion: if a concept (conception) of set also incorporates an account of classes, then the Intrinsicness Issue for large cardinals becomes quite a different issue. But, then, it is motivating classes \textit{intuitively} that becomes very hard. If, on the other hand, classes are just needed \textit{heuristically}, then one can always attempt to find equivalent non-class-involving formulations or, if one wishes to retain the vocabulary of classes \textit{as if} they were determinate entities, one can always do that with a view to using classes as a ladder, so to speak, to climb through sets (in particular, through large cardinals), and then get rid of the ladder at the end of the process. 

\subsection{Searching for Equilibrium}\label{equilibrium}

So far, we have considered the issue of whether large cardinals are \textit{intrinsically justified} taking for granted that the notion of `intrinsic justification' is fully understood. However, several authors have argued that this is not the case, that, for instance, the `derivability of an axiom from the concept of set' is a poorly, if at all, understood concept. In the previous subsection, we have reviewed instances of this attitude that question the identification of the concept of set with IC. But the issue is more general, and invocations of such seemingly explanatory notions of intrinsicness as `self-evidence', or `logicality' may be equally inadequate. 

Two alternative approaches consist in, on the one hand, viewing the concept of `being intrinsically justified' as too much a demand on an axiom, and, as a consequence, in contenting oneself with establishing an axiom's intrinsic \textit{plausibility} and, on the other, in taking \textit{extrinsic} justifications to be more cogent than \textit{intrinsic} ones.    

The latter approach is too extensively explored in the literature to be, even briefly, touched on here.\footnote{For an overview of this approach, see, in particular, \cite{maddy1997} and \cite{maddy2011}.} The former approach has been investigated by  \cite{parsons2008}, which introduces the notion of \textit{intrinsically plausible} (as opposed to \textit{intrinsically justified}) statements to refer to those set-theoretic statements which 

\begin{quote}

..are accepted without carrying the argument any further. They strike those who make them as true or evident, prior to any sense of how to construct an argument for them [...] This is not to say that there is no argument for the statement at all, but the possible arguments are of a more indirect character or have premises that are no more evident. (\cite{parsons2008}, p. 319)

\end{quote}

Parsons argues that the notion of `intrinsically justified axiom' seems to embody, and carry with itself, a reference to standards of \textit{conceptual evidence} which can hardly be met by most set-theoretic axioms, let alone, one would add, Large Cardinal Axioms. On the contrary, as Parsons explains, an intrinsically plausible axiom is one which, although not justifiable on purely logical grounds, yet is `perceived' as being true, that is, is taken to be true on some intuitive, pre-theoretic, yet fully rational, grounds.\footnote{Cf. \cite{parsons2008}, p. 320ff..}  

Now, any rational account of the plausibility of a set-theoretic axiom will also include, as one of the axiom's fundamental virtues, its ability to \textit{systematise} our knowledge, which would express itself in the `dialectical interplay' the axiom would foster between \textit{higher-level generalisations} and \textit{lower-level statements} in the context of a process that \cite{parsons1995}, based on \cite{rawls1971}, defines of \textit{reflective equilibrium}.\footnote{Cf. \cite{parsons1995}, pp. 69ff. Rawls' notion of reflective equilibrium is also discussed in great detail by \cite{hauser2001}, section 6 and \cite{hauser2002}, p. 275, in the context of Hauser's own conception of `theory formation'.}

So, for our purposes, one may argue that, if not intrinsically justified, Large Cardinal Axioms could, at least, be intrinsically plausible, in Parsons' sense, insofar as SRPs are intrinsically plausible. One of the fundamental reasons for this would, precisely, be the fact that SRPs have the patent ability to systematise the large-cardinal landscape to an unprecedented level of uniformity.  In particular, if Conjecture \ref{conj} is correct, then the study of large cardinals would, in practice, reduce to the study of SRPs. But there is another sense in which SRPs (hence, the large cardinals they yield) may be intrinsically plausible: they may also be an instance of higher-level generalisations interacting with lower-level set-theoretic (mathematical) statements, in two senses. 

In one sense, the Large Cardinal Axioms themselves, when formulated as  (\textit{local}) combinatorial statements may be taken to act as the relevant lower-level statements of which SRPs would represent the corresponding higher-level (\textit{global}) generalisations. On this picture, the meticulous matching of Large Cardinal Axioms with SRPs would represent a clear instance of reflective equilibrium in the sense specified.   

Alternatively, one could take any   mathematical statement (not necessarily related, in principle, to large cardinals) that follows from SRPs in a natural way as the lower-level statements interacting with SRPs (taken, again, as higher-order generalisations).  We are referring to mathematical results which emerge as more or less direct consequences of SRPs.

One immediate example is the fact, discussed in section \ref{just} and Remark \ref{LST}, that SR provides analogues (and strengthenings) of the Downwards-L\"owen-heim-Skolem Theorem. Other examples have emerged in areas of mathematics which seem remote from large cardinal theory. In \cite{bcmr2015}, for instance, the authors show that certain results in \textit{homotopy} and \textit{category} theory follow naturally from  SR restricted to definable classes of structures of an appropriate level of complexity ($\Pi_n$), which turns out to be equivalent to the existence of particular large cardinals (supercompact cardinals, for $n=1$, and $C^{n-1}$-extendible cardinals, for $n>1$) and which form a hierarchy, with VP now being equivalent to the existence of a $C^{n}$-extendible cardinal, for all $n$. Thus,  the hierarchy of such cardinals is equivalent to the hierarchy of  SRPs, and yield, in a uniform way, a hierarchy of results about, e.g.,  \emph{reflectivity} or \emph{smallness} of orthogonality classes for locally-presentable categories of structures, or \emph{cohomological localizations} of simplicial sets.\footnote{Further results along these lines on the existence of \textit{cohomological localisations} may be found in \cite{casacuberta2023}.}

\medskip

\cite{koellner2014} also takes Parsons' notion of intrinsically plausible (\textit{vis-à-vis} intrinsically justified) axiom as the starting point of a picture of justification hinged on what the author calls `evidentness order'. In essence, Koellner conjectures that, for each interpretability degree of theories $T$,\footnote{Two theories $T_1$ and $T_2$ are mutually interpretable if and only if, for all $\varphi$, if $T_1 \vdash \varphi$, then there is a translation $\tau(\varphi)$ in the language of $T_2$, such that $T_2 \vdash \tau(\varphi)$, and vice versa. The \textit{interpretability degree} of $T$ is the \textit{equivalence class} of all theories mutually interpretable with $T$. Cf. \cite{koellner2014}, p. 5.} one may find statements, say $\phi$ and $\psi$, such that, over $T$, $\phi$ is provably equivalent to $\psi$, but is \textit{less} evident than $\psi$, mostly in the sense, as explained by the author, that $\psi$ is more \textit{intrinsically plausible} than $\varphi$. 

One of Koellner's examples is, over some theory $T$ of the appropriate interpretability degree, the following pair of proof-theoretically equivalent statements: 1) (the Hydra Theorem): `Hercules will eventually chop off \textit{all} of the Hydra's heads after $n$ steps', and 2): `$\epsilon_0$ is well-ordered with respect to the primitive recursive relations', where 2) is, arguably, \textit{more evident} (more \textit{intrinsically plausible}) than 1). Then, one could say that, over $T$, 2) is `lower' than 1) in the evidentness order.\footnote{Cf. \cite{koellner2014}, pp. 12-16.} At this point, one may conjecture that, for each interpretability degree, there may be statements that are minimal in the evidentness order. Ultimately, such statements may be taken to be good approximations to \textit{intrinsically justified} statements.

Now, one could argue that SRPs are such statements. To begin with, one could argue that, for each appropriate interpretability degree, Large Cardinal Axioms are less evident (less intrinsically plausible) than the corresponding SRPs. The challenge would, then, be to show that, for each interpretability degree, there couldn't possibly be any statements lower than them in the evidentness order. In a sense, this would reduce to showing that there couldn't possibly be any statements equivalent to Large Cardinal Axioms \textit{more evident than} SRPs, in the sense explained in this section, and the work in the previous sections may be taken to have provided, at least, some evidence to this effect.

However, as far as we can see, there are two main problems with this strategy (and, as a consequence, potentially also with our strategy): first, it is not clear that a notion of `minimality', for each interpretability degree, may be sensibly formulated; second, it is not clear that the evidentness order is \textit{total}, which means that it is not necessarily the case that, for each interpretability degree, for any two statements $\phi$ and $\psi$, these are comparable; at least, no general methodology seems to be, so far, available to carry out this task in a principled way.\footnote{As acknowledged by Koellner himself in \cite{koellner2014}, p. 16, there might be pairs of statements for which there is, in fact, a lot of disagreement over which of the two is \textit{more evident}.}

\section{Concluding Remarks}

It might really happen that large cardinals will, eventually, be seen as intrinsically justified principles of set theory, possibly also as a result of the further considerations made in section \ref{equilibrium}. Until then, the case that they be so seems to be only partly conclusive. However, the good news is that the case is not entirely hopeless: in particular, in the paper we have shown that concrete progress can be made, through refining and improving on the existing mathematical principles based, in turn, on abstract motivating principles.    

Our main case study, throughout, has been Structural Reflection, in fact, \textit{the} Structural Reflection Principles. With reference to the Intrinsicness Issue, we have shown that, although SR, a clear instance of Resemblance, is firmly clung to the resources of IC, there are aspects of the mathematical formulation of SRPs which do not seem to automatically fall within the compass of IC, so SR's response to the issue is only partly satisfactory. In particular, as made clear in sections \ref{self} and \ref{just}, IC's sole resources do not seem to be sufficient to provide intrinsic motivation for \textit{elementary embeddings}, insofar as, essentially, IC is not able to provide intrinsic motivation for Self. However, as has been shown, one of the main hindrances to solving the issue, Schindler's Challenge, is successfully overcome by SR, so progress has, undeniably, been made.  

On the other hand, as far as the Universality Issue is concerned, SRPs represent a clear major breakthrough, insofar as they might be able to capture \textit{all} known large cardinal notions. As a consequence, they may, also, fill up an outstanding, and somewhat embarrassing, definitional void in the theory of large cardinals, i.e., the definition of `large cardinal' itself. 

It is, then, natural to foresee that the next task for the theory of Structural Reflection will consist in addressing the following question:

\begin{question}

Are Conjectures \ref{conj} and \ref{conj2} true? What (further) mathematical evidence would be needed to establish the truth of both?
    
\end{question}

On the purely philosophical side, what has emerged from the examination of the abstract motivating principles in section \ref{sec3} is that hardly any of these leads to the formulation of mathematical principles compatible with IC; in fact, all of these seem, in a way or the other, to essentially transcend the iterative conception and rather focus on the features of the universe of sets $V$ mostly taken to be a \textit{completed} totality. As a consequence, a philosophical (and, to some extent, also mathematical) question that one would think worth addressing in this respect is: 

\begin{question}
    Does there exist an abstract motivating principle for large cardinals which is licensed by, and whose mathematical expressions are all compatible with, IC?
\end{question}

One final word on the broader issue of justification in set theory that has been touched on many times. In the paper, we have identified a peculiar conceptual relationship between three elements of the process of justification itself:

\begin{enumerate}

    \item The concept of set (IC)

    \item Motivating principles

    \item Mathematical principles
    
\end{enumerate}

\noindent
whereby each of the three elements, fruitfully, and dynamically, feeds back into the other two.

Now, a clear novelty of our philosophical account has consisted in taking 2., the motivating principles (Reflection, Resemblance, etc.) to play a \textit{mediating} role between 1., the concept of set and 3., the mathematical principles themselves. In future work, one should investigate exactly what is the kind of `role' that they accomplish, and why this role is important: we have assumed throughout that they were, in a sense, just one further expression of the concept of set, but their nature seems to be far more complex. 

Finally, in section \ref{equilibrium}, the paper has also suggested that further progress in the justification of large cardinals could be attained through weakening the notion of `intrinsicness' itself. For instance, one could just attempt to show that large cardinals are \textit{natural} and \textit{indispensable} additions to the current axioms of set theory, or that they also meet further intra-theoretical desiderata, such as \textit{explanatoriness}, and even aesthetic ones such as \textit{beauty}.\footnote{For consideraions on these aspects of the large cardinal theory, see \cite{bagaria2005}, especially section 4.} But this would hardly count as justifying them intrinsically. Instead, the considerations made in section \ref{equilibrium}, based on Parsons' and Koellner's work on intrinsic plausibility, on the one hand, help reframe and reorient the problem of intrinsic justification in a (possibly) more productive way and, on the other, encourage set-theorists to find further lower-level statements directly entailed by SRPs that may help one reach a reflective equilibrium. These last considerations very naturally motivate the following further question:

\begin{question}
 To what extent could SRPs be further motivated by criteria and conceptions other than IC?   
\end{question}

\medskip

No doubt, there is still a lot to think about the intrinsic justifiability of large cardinals. The main merit of the paper, it seems to us, has consisted in showing that there exists a methodology that is already able to deliver concrete solutions and shed light on the fundamental philosophical issues involved. 

\pagebreak

\bibliography{Bib1}

\bibliographystyle{apalike}

\pagebreak

\tableofcontents

\end{document}